\documentclass[graybox]{svmult}

\usepackage{type1cm}
\usepackage{ dsfont }
\usepackage{makeidx}         
\usepackage{graphicx}        
\usepackage{multicol}        
\usepackage[bottom]{footmisc}

\usepackage{newtxtext}       %
\usepackage{newtxmath}

\usepackage{amsfonts}
\usepackage[ruled, vlined]{algorithm2e}
\usepackage[shortcuts]{extdash}

\def\ag#1{{\color{black}#1}}
\def\at#1{{\color{black}#1}}

 \makeindex

\title*{Unifying Framework for Accelerated Randomized Methods in Convex Optimization \thanks{Submitted to the editors {November 30, 2020}.
		{The research
		was supported by the Ministry of Science and Higher Education of the Russian Federation (Goszadaniye) No. 075-00337-20-03, project No. 0714-2020-0005.}}}

\author{Pavel Dvurechensky
	\thanks{Weierstrass Institute for Applied Analysis and Stochastics, Berlin; Moscow Institute of Physics and Technology, Moscow; HSE University, 
		Moscow, \email{pavel.dvurechensky@wias-berlin.de}}
	\and
	Alexander Gasnikov
	\thanks{Moscow Institute of Physics and Technology, Moscow; HSE University, Moscow; Weierstrass Institute for Applied Analysis and Stochastics, Berlin, 
		 \email{gasnikov.av@mipt.ru}}
	\and
	Alexander Tyurin
	\thanks{HSE University, Moscow, \email{alexandertiurin@gmail.com}}
	\and
	Vladimir Zholobov 
	\thanks{Moscow Institute of Physics and Technology, Moscow,
	\email{zholobov.va@phystech.edu}}
}

\newcommand{\R}{{\mathbb R}} 
\newcommand{\mE}{{\mathbb E}}
\renewcommand{\E}{{\mathbb E}}
\newcommand{\Sp}{{\mathcal S}}

\newcommand{\e}{\varepsilon}
\newcommand{\la}{\langle}
\newcommand{\ra}{\rangle}

\def\Rf{\mathcal R_p}
\def\Rb{\mathcal R_r}

\def\tf{\tilde{f}}

\def\tnf{\widetilde{\nabla}f}
\def\tnfg{\widetilde{\nabla}{(f+g)}}

\newcommand{\norm}[1]{\left\lVert#1\right\rVert}

\newcommand\uprule{\rule{0mm}{1.9ex}}
\newcommand{\argmin}{\operatornamewithlimits{argmin}}

\spnewtheorem{assumption}[definition]{Assumption}{}{}

\begin{document}

\maketitle

\abstract{
	In this paper, we consider smooth convex optimization problems with simple constraints and inexactness in the oracle information such as value, partial or directional derivatives of the objective function. We introduce a unifying framework, which allows to construct different types of accelerated randomized methods for such problems and to prove convergence rate theorems for them. We focus on accelerated random block-coordinate descent, accelerated random directional search, accelerated random derivative-free method and, using our framework, provide their versions for problems with inexact oracle information. Our contribution also includes accelerated random block-coordinate descent with inexact oracle and entropy proximal setup as well as derivative-free version of this method. Moreover, we present an extension of our framework for strongly convex optimization problems.
}


\section*{Introduction}
\label{S:Intro}
In this paper, we consider smooth convex optimization problems with simple constraints and inexactness in the oracle information such as objective value, partial or directional derivatives of the objective function.
Different types of randomized optimization algorithms, such as random coordinate descent or stochastic gradient descent for empirical risk minimization problem, have been extensively studied in the past decade with the main application being large-scale convex optimization problems. Our main focus in this paper is on accelerated randomized methods: random block-coordinate descent, random directional search, random derivative-free method. As opposed to non-accelerated methods, these methods have complexity proportional to $\frac{1}{\sqrt{\e}}$ iterations to achieve objective function residual $\e$, as opposed to $\frac{1}{\e}$ for non-accelerated methods. Accelerated random block-coordinate descent method was first proposed in \cite{nesterov2012efficiency}, which was the starting point for active research in this direction. The idea of the method is, on each iteration, to randomly choose a block of coordinates in the decision variable and make a step using the derivative of the objective function with respect to the chosen coordinates. Accelerated random directional search and accelerated random derivative-free method were first proposed in 2011 and published recently in \cite{nesterov2017random}, but there was no extensive research in this direction. The idea of random directional search is to use a projection of the objective's gradient onto a randomly chosen direction to make a step on each iteration. Random derivative-free method uses the same idea, but random projection of the gradient is approximated by finite-difference, i.e. the difference of values of the objective function in two close points. This also means that it is a zero-order method which uses only function values to make a step.  

Existing accelerated randomized methods have different convergence analysis. This motivated us to pose the main question, we address in this paper, as follows. \textit{Is it possible to find a crucial part of the convergence rate analysis and use it to systematically construct new accelerated randomized methods?} To some extent, our answer is "yes". We determine three main assumptions and use them to prove convergence rate theorem for our generic accelerated randomized method. Our framework allows both to reproduce known and to construct new accelerated randomized methods. The latter include new accelerated random block-coordinate descent with inexact block derivatives and entropy proximal setup.

\subsection*{Related Work}
In the seminal paper \cite{nesterov2012efficiency}, the author proposed random block-coordinate descent for convex optimization problems with simple convex separable constraints and accelerated random block-coordinate descent for unconstrained convex optimization problems. In \cite{lee2013efficient}, the authors proposed accelerated random block-coordinate descent with non-uniform probability of choosing a particular block of coordinates. They also developed an efficient implementation without full-dimensional operations on each iteration. The authors of \cite{fercoq2015accelerated} introduced accelerated block-coordinate descent for composite optimization problems, which include problems with separable constraints. Later, the paper \cite{lin2014accelerated} extended this method for strongly convex problems.  
The papers \cite{nesterov2017efficiency} \ag{(this work first appeared in May, 2015)}, \cite{allen2016even} and \cite{gasnikov2016accelerated} proposed an accelerated block-coordinate descent with complexity, which does not explicitly depend on the problem dimension. 
We also mention special type of accelerated block-coordinate descent developed in \cite{shalev-shwartz2014accelerated} for empirical risk minimization problems. All these accelerated block-coordinate descent methods work in Euclidean setup, when the norm in each block is Euclidean and is defined using some positive semidefinite matrix. Non-accelerated block-coordinate methods, but with non-euclidean setup, were considered in \cite{dang2015stochastic}. All the mentioned methods rely on exact block derivatives and exact projection on each step. Inexact projection in the context of non-accelerated random coordinate descent was considered in \cite{tappenden2016inexact}.

Research on accelerated random directional search and accelerated random derivative\-/free methods for smooth problems started in \cite{nesterov2017random}, where also non-smooth problems were considered. Derivative-free methods for non-smooth problems were further developed in \cite{agarwal2010optimal,duchi2015optimal,shamir2017optimal} for the case of exact computations and in \cite{gasnikov2016gradient-free, gasnikov2016stochastic,gasnikov2017stochastic,bayandina2017gradient-free} for the case of inexact function values. Non-accelerated gradient-free methods for smooth problems were developed in \cite{ghadimi2016mini-batch,ghadimi2013stochastic} for the case of exact objective value evaluation and in \cite{bogolubsky2016learning,berahas2019global,berahas2020theoretical} for inexact values, see also extensive reviews \cite{conn2009introduction,Larson_2019}. Accelerated gradient-free methods with inexact oracle were introduced in \cite{dvurechensky2018accelerated} and later extended for the case of random directional search in \cite{dvurechensky2020accelerated}. The main difference with this work is that here we develop a unifying framework for all the three types of accelerated randomized methods including coordinate descent method.

We should also mention that there are other accelerated randomized methods in \cite{frostig2015un-regularizing,lin2015universal,yuchen2015stochastic,allen2016katyusha,lan2017optimal}. Most of these methods were developed deliberately for empirical risk minimization problems and do not fall in the scope of this paper. We also list the following works, which consider accelerated full-gradient methods with inexact oracle \cite{aspremont2008smooth,devolder2014first,dvurechensky2016stochastic,cohen2018acceleration}.
%
%
%
%

\subsection*{Our Approach and Contributions}
Our framework has two main components, namely, Randomized Inexact Oracle and Unified Accelerated Randomized Method. The starting point for the definition of our oracle is a unified view on random directional search and random block-coordinate descent. In both these methods, on each iteration, a randomized approximation for the objective function's gradient is calculated and used, instead of the true gradient, to make a step. This approximation for the gradient is constructed by a projection on a randomly chosen subspace. For random directional search, this subspace is the line along a randomly generated direction. As a result a directional derivative in this direction is calculated. For random block-coordinate descent, this subspace is given by randomly chosen block of coordinates and block derivative is calculated. One of the key features of these approximations is that they are unbiased, i.e. their expectation is equal to the true gradient. We generalize two mentioned approaches by allowing other types of random transformations of the gradient for constructing its randomized approximation. 

The inexactness of our oracle is inspired by the relation between derivative-free method and directional search. In the framework of derivative-free methods, only the value of the objective function is available for use in an algorithm. At the same time, if the objective function is smooth, the directional derivative can be well approximated by the difference of function values at two points which are close to each other. Thus, in the context of zero-order optimization, one can calculate only an inexact directional derivative. Hence, one can construct only a biased randomized approximation for the gradient when a random direction is used. We combine previously mentioned random transformations of the gradient with possible inexactness of this transformations to construct our \textit{Randomized Inexact Oracle}, which we use in our generic algorithm to make a step on each iteration. 

The basis of our generic algorithm is Similar Triangles Method of \cite{tyurin2017mirror} (see also \cite{dvurechensky2020stable}), which is an accelerated gradient method with only one proximal mapping on each iteration, this proximal mapping being essentially the Mirror Descent step. The notable point is that, we only need to substitute the true gradient with our Randomized Inexact Oracle and slightly change one step in the Similar Triangles Method, to obtain our generic accelerated randomized algorithm, which we call Unified Accelerated Randomized Method (UARM), see \ag{Algorithm}~\ref{Alg:UARM}.
We prove convergence rate theorem for UARM in two cases: the inexactness of Randomized Inexact Oracle can be controlled and adjusted on each iteration of the algorithm, the inexactness can not be controlled.

We apply our framework to several particular settings: random directional search, random coordinate descent, random block-coordinate descent and their combinations with derivative-free approach. As a corollary of our main theorem, we obtain both known and new results on the convergence of different accelerated randomized methods with inexact oracle. 

To sum up, our contributions in this paper are as follows.
\begin{itemize}
	\item We introduce a general framework for constructing and analyzing different types of accelerated randomized methods, such as accelerated random directional search, accelerated block-coordinate descent, accelerated derivative-free methods. Our framework allows to obtain both known and new methods and their convergence rate guarantees as a corollary of our main \ag{Theorem}~\ref{Th:1}. 
	\item Using our framework, we introduce new accelerated methods with inexact oracle, namely, accelerated random directional search, accelerated random block-coordinate descent, accelerated derivative-free method. To the best of our knowledge, such methods with inexact oracle were not known before. See section~\ref{S:Appl}.
	\item Based on our framework, we introduce new accelerated random block\-/coordinate descent with inexact oracle and non-euclidean setup, which was not done before in the literature. The main application of this method is minimization of functions on a direct product of large number of low-dimensional simplexes. See subsection~\ref{S:RBCD}.
	\item We introduce new accelerated random derivative-free block-coordinate descent with inexact oracle and non-euclidean setup. Such method was not known before in the literature. Our method is similar to the method in the previous item, but uses only finite-difference approximations for block derivatives. See subsection~\ref{S:RDFBCD}.
\end{itemize}

The rest of the paper is organized as follows. In section~\ref{S:prel}, we provide the problem statement, motivate and make three our main assumptions, illustrate them by random directional search and random block-coordinate descent. In section~\ref{S:UARM}, we introduce our main algorithm, called Unified Accelerated Randomized Method, and, based on stated general assumptions, prove convergence rate \ag{in Theorem}~ \ref{Th:1}. Section~\ref{S:Appl} is devoted to applications of our general framework for different particular settings, namely
\begin{itemize}
	\item Accelerated Random Directional Search (\ref{S:RDS}),
	\item Accelerated Random Coordinate Descent (\ref{S:RCD}),
	\item Accelerated Random Block-Coordinate Descent (\ref{S:RBCD}),
	\item Accelerated Random Derivative-Free Directional Search (\ref{S:RDFDS}),
	\item Accelerated Random Derivative-Free Coordinate Descent (\ref{S:RDFCD}),
	\item Accelerated Random Derivative-Free Block-Coordinate Descent 
	(\ref{S:RDFBCD}).
	\item Accelerated Random Derivative-Free Block-Coordinate Descent with Random Approximations for Block Derivatives (\ref{S:RDFBCDRA}).
\end{itemize}

\section{Preliminaries}
\label{S:prel}

\subsection{Notation}
\label{S:Not}
Let finite-dimensional real vector space $E$ be a direct product of $n$ finite-dimensional real vector spaces $E_i$, $i=1,...,n$, i.e. $E=\otimes_{i=1}^n E_i$ and $\text{dim} E_i = p_i$, $i=1,...,n$. Denote also $p=\sum_{i=1}^n p_i$. Let, for $i=1,...,n$, $E_i^*$ denote the dual space for $E_i$. Then, the space dual to $E$ is $E^*=\otimes_{i=1}^n E_i^*$. Given a vector $x^{(i)} \in E_i$ for some $i \in 1,...,n$, we denote as $[x^{(i)}]_j$ its $j$-th coordinate, where $j \in 1,...,p_i$.
To formalize the relationship between vectors in $E_i$, $i=1,...,n$ and vectors in $E$, we define primal partition operators $U_i : E_i \to E$, $i=1,...,n$, by identity
\begin{equation}
x = (x^{(1)},...,x^{(n)}) = \sum_{i=1}^n U_i x^{(i)}, \quad x^{(i)} \in E_i, \quad i=1,...,n, \quad x \in E.
\label{eq:UDef}
\end{equation}
For any fixed $i \in 1,...,n$, $U_i$ maps a vector $x^{(i)} \in E_i$, to the vector $(0,....,x^{(i)},...,0) \in E$. The adjoint operator $U_i^T: E^* \to E_i^*$, then, is an operator, which, maps a vector $g = (g^{(1)},...,g^{(i)},..., g^{(n)}) \in E^*$, to the vector $g^{(i)} \in E_i^*$. 
Similarly, we define dual partition operators $\widetilde{U}_i : E_i^* \to E^*$, $i=1,...,n$, by identity
\begin{equation}
g = (g^{(1)},...,g^{(n)}) = \sum_{i=1}^n \widetilde{U_i} g^{(i)}, \quad g^{(i)} \in E_i^*, \quad i=1,...,n, \quad g \in E^*.
\label{eq:tUDef}
\end{equation}
For all $i=1,...,n$, we denote the value of a linear function $g^{(i)} \in E_i^*$ at a point $x^{(i)}\in E_i$ by $\la g^{(i)}, x^{(i)} \ra_i$. We define 
$$
\la g, x \ra = \sum_{i=1}^n \la g^{(i)}, x^{(i)} \ra_i, \quad x \in E, \quad g \in E^*.
$$
For all $i=1,...,n$, let $\|\cdot\|_i$ be some norm on $E_i$ and $\|\cdot\|_{i,*}$ be the norm on $E_i^*$ which is dual to $\|\cdot\|_i$
$$
\|g^{(i)}\|_{i,*} = \max_{\|x^{(i)}\|_i \leq 1} \la g^{(i)}, x^{(i)} \ra_i.
$$
Given parameters $\beta_i \in \R^n_{++}$, $i=1,...,n$, we define the norm of a vector $x=(x^{(1)},...,x^{(n)}) \in E$ as 
$$
\|x\|_E^2 = \sum_{i=1}^n \beta_i \|x^{(i)}\|_i^2.
$$ 
Then, clearly, the dual norm of a vector $g=(g^{(1)},...,g^{(n)}) \in E^*$ is 
$$
\|g\|_{E,*}^2 = \sum_{i=1}^n \beta_i^{-1} \|g^{(i)}\|_i^2.
$$ 
Throughout the paper, we consider optimization problem with feasible set $Q$, which is assumed to be given as $Q=\otimes_{i=1}^n Q_i \subseteq E$, where $Q_i \subseteq E_i$, $i=1,...,n$ are closed convex sets.  
To have more flexibility and be able to adapt algorithm to the structure of sets $Q_i$, $i=1,...,n$, we introduce {\it proximal setup}, see e.g. \cite{ben-tal2015lectures}.
For all $i=1,...,n$, we choose a {\it prox-function} $d_i(x^{(i)})$ which is continuous, convex on $Q_i$ and
\begin{enumerate}
	\item admits a continuous in $x^{(i)} \in Q_i^0$ selection of subgradients 	$\nabla d_i(x^{(i)})$, where $x^{(i)} \in Q_i^0 \subseteq Q_i$, and $Q_i^0$  is the set of all $x^{(i)}$, where $\nabla d_i(x^{(i)})$ exists;
	\item is $1$-strongly convex on $Q_i$ with respect to $\|\cdot\|_i$, i.e., for any $x^{(i)} \in Q_i^0, y^{(i)} \in Q_i$, it holds that $d_i(y^{(i)})-d_i(x^{(i)}) -\la \nabla d_i(x^{(i)}) ,y^{(i)}-x^{(i)} \ra_i \geq \frac12\|y^{(i)}-x^{(i)}\|_i^2$.
\end{enumerate} 
We define also the corresponding {\it Bregman divergence} $V_i[z^{(i)}] (x^{(i)}) := d_i(x^{(i)}) - d_i(z^{(i)}) - \la \nabla d_i(z^{(i)}), x^{(i)} - z^{(i)} \ra_i$, $x^{(i)} \in Q_i, z^{(i)} \in Q_i^0$, $i=1,...,n$. 
It is easy to see that $V_i[z^{(i)}] (x^{(i)}) \geq \frac12\|x^{(i)} - z^{(i)}\|_i^2, \quad x^{(i)} \in Q_i, z^{(i)} \in Q_i^0$, $i=1,...,n$.
Standard proximal setups, e.g. Euclidean, entropy, $\ell_1/\ell_2$, simplex can be found in \cite{ben-tal2015lectures}. It is easy to check that, for given parameters $\beta_i \in \R^n_{++}$, $i=1,...,n$, the functions $d(x) = \sum_{i=1}^n \beta_i d_i(x^{(i)})$ and $V[z] (x) = \sum_{i=1}^n \beta_i V_i[z^{(i)}] (x^{(i)}) $ are respectively a prox-function and a Bregman divergence corresponding to $Q$. 
Also, clearly,
\begin{equation}
V[z] (x) \geq \frac12\|x - z\|_E^2, \quad x \in Q, \quad z \in Q^0 := \otimes_{i=1}^n Q_i^0.
\label{eq:BFLowBound}
\end{equation}
For a differentiable function $f(x)$, we denote by $\nabla f(x) \in E^*$ its gradient. 

\subsection{Problem Statement and Assumptions}
\label{S:PrSt&Asmpt}

The main problem, we consider, is as follows
\begin{equation}
\min_{x\in Q \subseteq E}  f(x), 
\label{eq:PrSt}
\end{equation}
where $f(x)$ is a smooth convex function, $Q=\otimes_{i=1}^n Q_i \subseteq E$, with $Q_i \subseteq E_i$, $i=1,...,n$ being closed convex sets. 

We now list our main assumptions and illustrate them by two simple examples. More detailed examples are given in section~\ref{S:Appl}. As the first example here, we consider random directional search, in which the gradient of the function $f$ is approximated by a vector $\la \nabla f(x) , e \ra e$, where $\la \nabla f(x) , e \ra$ is the directional derivative in direction $e$ and random vector $e$ is uniformly distributed over the Euclidean sphere of radius 1. Our second example is random block-coordinate descent, in which the gradient of the function $f$ is approximated by a vector $\widetilde{U}_iU_i^T \nabla f(x)$, where $U_i^T \nabla f(x)$ is $i$-th block derivative and the block number $i$ is uniformly randomly sampled from $1,...,n$. The common part in both these randomized gradient approximations is that, first, one randomly chooses a subspace which is either the line, parallel to $e$, or $i$-th block of coordinates. Then, one projects the gradient on this subspace by calculating either $\la \nabla f(x) , e \ra $ or $U_i^T \nabla f(x)$. Finally, one lifts the obtained random projection back to the whole space $E$ either by multiplying directional derivative by vector $e$, or applying dual partition operator $\widetilde{U}_i$. At the same time, in both cases, if one scales the obtained randomized approximation for the gradient by multiplying it by $n$, one obtains an \textit{unbiased randomized approximation} of the gradient
$$
\E_e n \la \nabla f(x) , e \ra e = \nabla f(x), \quad \mE_i n \widetilde{U}_iU_i^T \nabla f(x) = \nabla f(x), \quad x \in Q.
$$
We also want our approach to allow construction of derivative-free methods. For a function $f$ with $L$-Lipschitz-continuous gradient, the directional derivative can be well approximated by the difference of function values in two close points. Namely, it holds that 
$$
\la \nabla f(x) , e \ra  = \frac{f(x+\tau e) - f(x)}{\tau} + o(\tau),
$$
where $\tau > 0$ is a small parameter. Thus, if only the value of the function is available, one can calculate only inexact directional derivative, which leads to biased randomized approximation for the gradient if the direction is chosen randomly. These three features, namely, random projection and lifting up, unbiased part of the randomized approximation for the gradient, bias in the randomized approximation for the gradient, lead us to the following assumption about the structure of our general \textit{Randomized Inexact Oracle}.

\begin{assumption}[Randomized Inexact Oracle]\label{A:RIO}
	
	We access the function $f$ only through Randomized Inexact Oracle $\tnf(x)$, $x \in Q$, which is given by
	\begin{equation}
	\tnf(x) = \rho \Rb(\Rf^T\nabla f(x) + \xi(x)) \in E^*,
	\label{eq:IODef}
	\end{equation}
	where $\rho> 0$ is a known constant; $\Rf$ is a random "`projection"' operator from some auxiliary space $H$ to $E$, and, hence, $\Rf^T$, acting from $E^*$ to $H^*$, is the adjoint to $\Rf$; $\Rb: H^* \to E^*$ is also some random "`reconstruction"' operator; $\xi(x) \in H^*$ is a, possibly random, vector characterizing the error of the oracle. The oracle is also assumed to satisfy the following properties
	\begin{align}
	&\mE \rho \Rb\Rf^T\nabla f(x) = \nabla f(x), \quad \forall x \in Q, 
	\label{eq:IOProp1} \\
	&\|\Rb \xi(x)\|_{E,*} \leq \delta, \quad \forall x \in Q,
	\label{eq:IOProp2}
	\end{align} 
	where $\delta \geq 0$ is oracle \textit{error level}.
\end{assumption}

Let us make some comments on this assumption. The nature of the operator $\Rf^T$ is generalization of random projection. 
For the case of random directional search, $H = \R$, $\Rf^T: E^* \to \R$ is given by $\Rf^T g = \la g, e \ra$, $g \in E^*$. For the case of random block-coordinate descent, $H = E_i$, $\Rf^T: E^* \to E_i^*$ is given by $\Rf^T g = U_i^T g$, $g \in E^*$. We assume that there is some additive error $\xi(x)$ in the generalized random projection $\Rf^T \nabla f(x)$. This error can be introduced, for example, when finite-difference approximation of the directional derivative is used. Finally, we lift the inexact random projection $\Rf^T \nabla f(x)+\xi(x)$ back to $E$ by applying operator $\Rb$. For the case of random directional search, $\Rb: \R \to E^*$ is given by $\Rb t = te$, $t \in \R$. For the case of random block-coordinate descent, $\Rb: E_i^* \to E^*$ is given by $\Rb g^{(i)} = \widetilde{U}_i g^{(i)}$, $g^{(i)} \in E_i^*$. The number $\rho$ is the normalizing coefficient, which allows the part $\rho \Rb\Rf^T\nabla f(x)$ to be unbiased randomized approximation for the gradient. This is expressed by equality ~\eqref{eq:IOProp1}. Finally, we assume that the error in our oracle is bounded, which is expressed by property ~\eqref{eq:IOProp2}. In our analysis, we consider two cases: the error $\xi$ can be controlled and $\delta$ can be appropriately chosen on each iteration of the algorithm; the error $\xi$ can not be controlled and we only know oracle error level $\delta$.

Let us move to the next assumption. As said, our generic algorithm is based on Similar Triangles Method of \cite{gasnikov2018universal_,tyurin2017mirror} (see also \cite{dvurechensky2020stable}), which is an accelerated gradient method with only one proximal mapping on each iteration. This proximal mapping is essentially the Mirror Descent step. For simplicity, let us consider here an unconstrained minimization problem in the Euclidean setting. This means that $Q_i = E_i = \R^{p_i}$, $\|x^{(i)}\|_i = \|x^{(i)}\|_2$, $i=1,...,n$. Then, given a point $u \in E$, a number $\alpha$, and the gradient $\nabla f(y)$ at some point $y \in E$, the Mirror Descent step is
$$
u_+ = \arg \min_{x \in E} \left\{\frac{1}{2}\|x-u\|_2^2 + \alpha \la \nabla f(y),x \ra\right\} = u - \alpha \nabla f(y).
$$
Now we want to substitute the gradient $\nabla f(y)$ with our Randomized Inexact Oracle $\tnf(y)$. Then, we see that the step $u_+ = u - \alpha \tnf(y)$ makes progress only in the subspace onto which the gradient is projected, while constructing the Randomized Inexact Oracle. In other words, $u - u_+$ lies in the same subspace as $\tnf(y)$. In our analysis, this is a desirable property and we formalize it as follows.
\begin{assumption}[Regularity of Prox-Mapping]
	\label{A:RPM}
	The set $Q$, norm $\|\cdot\|_E$, prox-function $d(x)$, and Randomized Inexact Oracle $\tnf(x)$ are chosen in such a way that, for any $u,y \in Q$, $\alpha >0$, the point
	\begin{equation}
	u_+ = \arg \min_{x \in Q} \left\{V[u](x) + \alpha \la \tnf(y),x \ra\right\}
	\label{eq:uPlusDef}
	\end{equation}
	satisfies 
	\begin{equation}
	\la \Rb\Rf^T\nabla f(y), u - u_+ \ra = \la \nabla f(y), u - u_+ \ra.
	\label{eq:RPMA}
	\end{equation}
\end{assumption}
The interpretation is that, in terms of linear pairing with $u - u_+$, the unbiased part $\Rb\Rf^T\nabla f(y)$ of the Randomized Inexact Oracle makes the same progress as the true gradient $\nabla f(y)$.

Finally, we want to formalize the smoothness assumption for the function $f$. In our analysis, we use only the smoothness of $f$ in the direction of $u_+ - u$, where $u \in Q$ and $u_+$ is defined in ~\eqref{eq:uPlusDef}. Thus, we consider two points $x,y \in Q$, which satisfy equality $x = y + a (u_+ - u)$, where $a \in \R$. For the random directional search, it is natural to assume that $f$ has $L$-Lipschitz-continuous gradient with respect to the Euclidean norm, i.e. 
\begin{equation}
f(x) \leq f(y) + \la \nabla f(y), x-y \ra + \frac{L}{2}\|x-y\|_2^2, \quad x,y \in Q.
\label{eq:fLipSm0}
\end{equation}
Then, if we define $\|x\|_E^2 = L \|x\|_2^2$, we obtain that, for our choice $x = y + a (u_+ - u)$,
$$
f(x) \leq f(y) + \la \nabla f(y), x - y \ra + \frac12\|x-y\|_E^2.
$$
Usual assumption for random block-coordinate descent is that the gradient of $f$ is block-wise Lipschitz continuous. This means that, for all $i=1,...,n$, block derivative $f'_i(x) = U_i^T \nabla f(x)$ is $L_i$-Lipschitz continuous with respect to chosen norm $\|\cdot \|_i$, i.e.
\begin{equation}
\|f'_i(x+U_ih^{(i)}) - f'_i(x)\|_{i,*} \leq L_i \|h^{(i)}\|_i, \quad h^{(i)} \in E_i, \quad i=1,...,n, \quad x \in Q.
\label{eq:fBlockLipSm0}
\end{equation}
By the standard reasoning, using ~\eqref{eq:fBlockLipSm0}, one can prove that, for all $i=1,...,n$,
\begin{equation}
f(x+U_ih^{(i)}) \leq f(x) + \la U_i^T\nabla f(x), h^{(i)} \ra + \frac{L_i}{2}\|h^{(i)}\|_i^2, \quad h^{(i)} \in E_i, \quad x \in Q.
\label{eq:fBlockLipSm00}
\end{equation}
In block-coordinate setting, $\tnf(x)$ has non-zero elements only in one, say $i$-th, block and it follows from \eqref{eq:uPlusDef} that $u_+ - u$ also has non-zero components only in the $i$-th block. Hence, there exists $h^{(i)} \in E_i$, such that $u_+ - u = U_i h_i $ and $x = y + a U_i h^{(i)} $.
Then, if we define $\|x\|_E^2 = \sum_{i=1}^n L_i \|x^{(i)}\|_i^2$, we obtain
\begin{align}
f(x)  = f(y + a U_i h^{(i)}) & \stackrel{~\eqref{eq:fBlockLipSm00}}{\leq}  f(y) + \la U_i^T\nabla f(y), a h^{(i)} \ra + \frac{L_i}{2}\|a h^{(i)}\|_i^2 \notag \\
& = f(y) + \la \nabla f(y), a U_i h^{(i)} \ra + \frac{1}{2}\|a U_i h^{(i)}\|_E^2  \notag \\
& = f(y) + \la \nabla f(y), x - y \ra + \frac{1}{2}\|x-y\|_E^2. \notag 
\end{align}
We generalize these two examples and assume smoothness of $f$ in the following sense.
\begin{assumption}[Smoothness]	\label{A:S}
    The norm $\|\cdot\|_E$ is chosen in such a way that, for any $u,y \in Q$, $a \in \R$, if $x = y + a (u_+ - u) \in Q$, then
	\begin{equation}
	f(x) \leq f(y) + \la \nabla f(y), x - y \ra + \frac12\|x-y\|_E^2.
	\label{eq:NfLipAsm}
	\end{equation}
\end{assumption}

\section{Unified Accelerated Randomized Method}
\label{S:UARM}
In this section, we introduce our generic Unified Accelerated Randomized Method, which is listed as algorithm~\ref{Alg:UARM} below, and prove \ref{Th:1}, which gives its convergence rate. The method is constructed by a modification of Similar Triangles Method (see \cite{dvurechensky2020stable}) and, thus, inherits part of its name. 
\begin{algorithm}[h!]
 \caption{Unified Accelerated Randomized Method (UARM)}
 \label{Alg:UARM}
 \SetKwInOut{Input}{Input}\SetKwInOut{Output}{Output}
 \KwIn{starting point $u_0 \in Q^0 = \otimes_{i=1}^n Q_i^0$, \at{number of iterations $N$,} prox-setup:  $d(x)$, $V[u] (x)$, see subsection~\ref{S:Not}.}
 Set $k=0$, $A_0=\alpha_0=1-\frac{1}{\rho}$, $x_0=y_0=u_0$.
 
 \Repeat{\ag{$k\le N$}}{
  Find $\alpha_{k+1}$ as the largest root of the equation
			\begin{equation}
			A_{k+1}:=A_k+\alpha_{k+1} = \rho^2 \alpha_{k+1}^2.
			\label{eq:alpQuadEq}
			\end{equation}
  Calculate 
			\begin{equation}
			y_{k+1} = \frac{\alpha_{k+1}u_k + A_k x_k}{A_{k+1}} .
			\label{eq:ykp1Def}
			\end{equation}
  Calculate 
			\begin{equation}
			u_{k+1}=\arg \min_{x \in Q} \{V[u_{k}](x) + \alpha_{k+1}\la \tnf(y_{k+1}), x \ra \}.
			\label{eq:ukp1Def}
			\end{equation}
 Calculate
			\begin{equation}
			x_{k+1} = y_{k+1} + \rho \frac{\alpha_{k+1}}{A_{k+1}}(u_{k+1} -u_k).
			\label{eq:xkp1Def}
			\end{equation}
 Set $k=k+1$.}
 \KwOut{The point \ag{$x_{N}$}.}
\end{algorithm}


\begin{lemma}
	Algorithm~\ref{Alg:UARM} is correctly defined in the sense that, for all $k\geq 0$, $x_k,y_k \in Q$.
	\label{Lm:ConvComb}
\end{lemma}

\begin{proof}
	The proof is a direct generalization of Lemma 2 in \cite{fercoq2015accelerated}. By definition ~\eqref{eq:ukp1Def}, for all $k\geq 0$, $u_k \in Q$. If we prove that, for all $k\geq 0$, $x_k \in Q$, then, from ~\eqref{eq:ykp1Def}, it follows that, for all $k\geq 0$, $y_k \in Q$. Let us prove that, for all $k\geq 0$, $x_k$ is a convex combination of $u_0 \dots u_{k}$, namely $x_{k} = \sum_{l=0}^{k}\gamma_{k}^l u_{l}$, where $\gamma_0^0 = 1$, $\gamma_1^0 = 0$, $\gamma_1^1 = 1$, and for $k \geq 1$,
	\begin{equation}
	\label{eq:gammakl}
	\gamma_{k+1}^l =\begin{cases} \left(1 - \frac{\alpha_{k+1}}{A_{k+1}}\right)\gamma_k^l, & l = 0,\dots,k - 1 \\ \frac{\alpha_{k+1}}{A_{k+1}} \left(1 - \rho \frac{\alpha_{k}}{A_{k}}\right) + \rho \left(\frac{\alpha_{k}}{A_{k}} - \frac{\alpha_{k+1}}{A_{k+1}}\right),& l = k \\
	\rho \frac{\alpha_{k+1}}{A_{k+1}},& l = k+1. \end{cases}
	\end{equation}
	Since, $x_0=u_0$, we have that $\gamma_0^0 = 1$. Next, by ~\eqref{eq:xkp1Def}, we have $x_1 = y_{1} + \rho\frac{\alpha_{1}}{A_{1}} (u_{1} - u_0) = u_0 + \rho\frac{\alpha_{1}}{A_{1}} (u_{1} - u_0) = (1 - \rho\frac{\alpha_{1}}{A_{1}}) u_0 + \rho\frac{\alpha_{1}}{A_{1}} u_{1}$. Solving the equation ~\eqref{eq:alpQuadEq} for $k=0$, and using the choice $\alpha_0 = 1 - \frac{1}{\rho}$, we obtain that $\alpha_1 = \frac{1}{\rho}$ and
	\begin{equation}
	\frac{\alpha_1}{A_1} \stackrel{~\eqref{eq:alpQuadEq}}{=}  \frac{\alpha_1}{\rho^2 \alpha_1^2} = \frac{1}{\rho}.
	\label{eq:al1divA1}
	\end{equation}
	Hence, $x_1=u_1$ and	$\gamma_1^0 = 0$, $\gamma_1^1 = 1$. Let us now assume that $x_{k} = \sum_{l=0}^{k}\gamma_{k}^l u_{l}$ and prove that $x_{k+1}$ is also a convex combination with coefficients, given by ~\eqref{eq:gammakl}.
	From ~\eqref{eq:ykp1Def}, ~\eqref{eq:xkp1Def}, we have 
	\begin{gather*}
	x_{k+1} = y_{k+1} + \rho\frac{\alpha_{k+1}}{A_{k+1}} (u_{k+1} - u_k) = \frac{\alpha_{k+1}u_k + A_k x_k}{A_{k+1}} + \rho\frac{\alpha_{k+1}}{A_{k+1}} (u_{k+1} - u_k) \\ =
	\frac{A_k}{A_{k+1}}x_k + \left(\frac{\alpha_{k+1}}{A_{k+1}} - \rho\frac{\alpha_{k+1}}{A_{k+1}}\right)u_k + \rho\frac{\alpha_{k+1}}{A_{k+1}}u_{k+1} \\ = \left(1 - \frac{\alpha_{k+1}}{A_{k+1}}\right) \sum_{l = 0}^{k} \gamma_k^l u_l + \left(\frac{\alpha_{k+1}}{A_{k+1}} - \rho\frac{\alpha_{k+1}}{A_{k+1}}\right)u_k + \rho\frac{\alpha_{k+1}}{A_{k+1}}u_{k+1}.
	\end{gather*}
	Note that all the coefficients sum to 1.
	Next, we have 
	\begin{align*}
	x_{k+1} =& 
	\left(1 - \frac{\alpha_{k+1}}{A_{k+1}}\right) \sum_{l = 0}^{k-1} \gamma_k^l u_l\\ &+ \left(\gamma_k^k \left(1 - \frac{\alpha_{k+1}}{A_{k+1}}\right) + \left(\frac{\alpha_{k+1}}{A_{k+1}} - \rho \frac{\alpha_{k+1}}{A_{k+1}}\right)\right)u_k + \rho \frac{\alpha_{k+1}}{A_{k+1}}u_{k+1}  \\
	&\hspace{-1em}=
	\left(1 - \frac{\alpha_{k+1}}{A_{k+1}}\right) \sum_{l = 0}^{k-1} \gamma_k^l u_l\\
	&+ \left(\rho  \frac{\alpha_{k}}{A_{k}} \left(1 - \frac{\alpha_{k+1}}{A_{k+1}}\right) + \left(\frac{\alpha_{k+1}}{A_{k+1}} - \rho \frac{\alpha_{k+1}}{A_{k+1}}\right)\right)u_k + \rho \frac{\alpha_{k+1}}{A_{k+1}}u_{k+1}  \\
	&\hspace{-1em}=\left(1 - \frac{\alpha_{k+1}}{A_{k+1}}\right) \sum_{l = 0}^{k-1} \gamma_k^l u_l\\ 
	&+ \left(\frac{\alpha_{k+1}}{A_{k+1}} \left(1 - \rho \frac{\alpha_{k}}{A_{k}}\right) + \rho \left(\frac{\alpha_{k}}{A_{k}} -\frac{\alpha_{k+1}}{A_{k+1}}\right)\right)u_k + \rho \frac{\alpha_{k+1}}{A_{k+1}}u_{k+1}.
	\end{align*}
	So, we see that ~\eqref{eq:gammakl} holds for $k+1$. It remains to show that $\gamma_{k+1}^l \geq 0$, $l = 0,\dots,k+1$.
	For $\gamma_{k+1}^l$, $l = 0,\ldots,k-1$ and $\gamma_{k+1}^{k+1}$ it is obvious. From ~\eqref{eq:alpQuadEq}, we have
	$$
	\alpha_{k+1} = \frac{1+\sqrt{1+4 \rho^2 A_k}}{2 \rho^2}.
	$$
	Thus, since $\{A_k\}$, $k\geq 0$ is non-decreasing sequence, $\{\alpha_{k+1}\}$, $k\geq 0$ is also non-decreasing. From ~\eqref{eq:alpQuadEq}, we obtain $\frac{\alpha_{k+1}}{A_{k+1}} = \frac{\alpha_{k+1}}{\rho^2\alpha_{k+1}^2}$, which means that this sequence is non-increasing. Thus, $\frac{\alpha_{k}}{A_{k}} \geq \frac{\alpha_{k+1}}{A_{k+1}}$ and $\frac{\alpha_{k}}{A_{k}} \leq \frac{\alpha_{1}}{A_{1}} \leq \frac{1}{\rho}$ for $k \geq 1$. These inequalities prove that $\gamma_{k+1}^{k} \geq 0$.	
\end{proof}

\begin{lemma}
	\label{Lm:1}
	Let the sequences $\{x_k, y_k, u_k, \alpha_k, A_k \}$, $k\geq 0$ be generated by algorithm~\ref{Alg:UARM}. Then,  for all $u \in Q$, it holds that
	\begin{align}
	\alpha_{k+1}\la \tnf(y_{k+1}), u_k - u\ra \leq & \,\,A_{k+1}(f(y_{k+1}) - f(x_{k+1})) + V[u_k](u) - V[u_{k+1}](u) \notag \\
	& + \alpha_{k+1} \rho \la \Rb\xi(y_{k+1}), u_k - u_{k+1} \ra.
	\label{eq:Lm1}
	\end{align}	
\end{lemma}

\begin{proof}
	Using assumptions~\ref{A:RIO} and \ref{A:RPM} with $\alpha = \alpha_{k+1}$, $y=y_{k+1}$, $u=u_k$, $u_+ = u_{k+1}$, we obtain
	\begin{align}
	\alpha_{k+1}\la \tnf(y_{k+1}), u_k - u_{k+1}\ra & \stackrel{~\eqref{eq:IODef}}{=} \alpha_{k+1} \rho \la \Rb(\Rf^T\nabla f(y_{k+1}) + \xi(y_{k+1})), u_k - u_{k+1} \ra \notag \\
	& \hspace{-5em}\stackrel{~\eqref{eq:RPMA}}{=} \alpha_{k+1} \rho \la \nabla f(y_{k+1}), u_k - u_{k+1} \ra  + \alpha_{k+1} \rho \la \Rb\xi(y_{k+1}), u_k - u_{k+1} \ra \notag \\
	& \hspace{-5em}\stackrel{~\eqref{eq:xkp1Def}}{=} A_{k+1} \la \nabla f(y_{k+1}), y_{k+1} - x_{k+1} \ra  + \alpha_{k+1} \rho \la \Rb\xi(y_{k+1}), u_k - u_{k+1} \ra  \label{eq:Lm1Pr1}.
	\end{align}	
	Note that, from the optimality condition in ~\eqref{eq:ukp1Def}, for any $u \in Q$, we have
	\begin{equation}
	\la \nabla V[u_{k}](u_{k+1}) + \alpha_{k+1} \tnf(y_{k+1}), u - u_{k+1} \ra \geq 0.
	\label{eq:Lm1Pr2}
	\end{equation}
	By the definition of $V[u](x)$, we obtain, for any $u \in Q$,
	\begin{align}
	V[u_k](u) - V[u_{k+1}](u) - V[u_k](u_{k+1})  = &d(u) - d(u_k) - \la \nabla d(u_k), u-u_k \ra \notag \\
	&- \left( d(u) - d(u_{k+1}) - \la \nabla d(u_{k+1}), u-u_{k+1} \ra  \right) \notag \\
	&- \left( d(u_{k+1}) - d(u_k) - \la \nabla d(u_k), u_{k+1}-u_k \ra \right) \notag \\
	& \hspace{-1em} = \la \nabla d(u_k) - \nabla d(u_{k+1}) , u_{k+1} - u \ra \notag \\
	& \hspace{-1em} = \la - \nabla V[u_{k}](u_{k+1}) , u_{k+1} - u \ra. 
	\label{eq:magicV}
	\end{align}
	Further, for any $u \in Q$, by assumption~\ref{A:S},
	\begin{align}
	\alpha_{k+1}\la \tnf(y_{k+1}), u_k - u\ra
	& = \alpha_{k+1}\la \tnf(y_{k+1}), u_k - u_{k+1}\ra + \alpha_{k+1}\la \tnf(y_{k+1}), u_{k+1} - u\ra \notag \\
	&\hspace{-7em} \stackrel{~\eqref{eq:Lm1Pr2}}{\leq} \alpha_{k+1}\la \tnf(y_{k+1}), u_k - u_{k+1}\ra + \la -\nabla V[u_{k}](u_{k+1}) , u_{k+1} - u  \ra \notag \\
	&\hspace{-7em} \stackrel{~\eqref{eq:magicV}}{=} \alpha_{k+1}\la \tnf(y_{k+1}), u_k - u_{k+1}\ra + V[u_k](u) - V[u_{k+1}](u) - V[u_k](u_{k+1}) \notag \\
	&\hspace{-7em} \stackrel{~\eqref{eq:BFLowBound}}{\leq} \alpha_{k+1}\la \tnf(y_{k+1}), u_k - u_{k+1}\ra + V[u_k](u) - V[u_{k+1}](u) - \frac12\|u_k-u_{k+1}\|_E^2 \notag \\\
	&\hspace{-8em} \stackrel{~\eqref{eq:Lm1Pr1},~\eqref{eq:xkp1Def}}{=} A_{k+1} \la \nabla f(y_{k+1}), y_{k+1} - x_{k+1} \ra  + \alpha_{k+1} \rho \la \Rb\xi(y_{k+1}), u_k - u_{k+1} \ra + \notag \\
	&\hspace{-5em} + V[u_k](u) - V[u_{k+1}](u) - \frac{A_{k+1}^2}{2\rho^2 \alpha_{k+1}^2}\|y_{k+1} - x_{k+1}\|_E^2 \notag \\
	&\hspace{-7em} \stackrel{~\eqref{eq:alpQuadEq}}{=} A_{k+1}\left(\la \nabla f(y_{k+1}), y_{k+1} - x_{k+1} \ra - \frac12\|y_{k+1} - x_{k+1}\|_E^2  \right) + \notag \\
	&\hspace{-5em} + \alpha_{k+1} \rho \la \Rb\xi(y_{k+1}), u_k - u_{k+1} \ra + V[u_k](u) - V[u_{k+1}](u) \notag \\
	&\hspace{-8em} \stackrel{~\eqref{eq:xkp1Def},~\eqref{eq:NfLipAsm}}{\leq} A_{k+1}\left(  f(y_{k+1}) - f(x_{k+1})  \right)  + V[u_k](u) - V[u_{k+1}](u) +\notag \\
	&\hspace{-5em} +\alpha_{k+1} \rho \la \Rb\xi(y_{k+1}), u_k - u_{k+1} \ra. \notag
	\end{align}
	In the last inequality, we used assumption~\ref{A:S} with $a=\rho \frac{\alpha_{k+1}}{A_{k+1}}$, $x=x_{k+1}$, $y=y_{k+1}$, $u=u_k$, $u_+=u_{k+1}$.
\end{proof}

\begin{lemma}
	\label{Lm:2}
	Let the sequences $\{x_k, y_k, u_k, \alpha_k, A_k \}$, $k\geq 0$ be generated by \ag{A}lgorithm~\ref{Alg:UARM}. Then,  for all $u \in Q$, it holds that
	\begin{align}
	\alpha_{k+1}\la \nabla f(y_{k+1}), u_k - u\ra \leq & \,\, A_{k+1}(f(y_{k+1}) - \E_{k+1}f(x_{k+1})) + V[u_k](u) \notag \\
	& - \E_{k+1}V[u_{k+1}](u) 
	+ 	\E_{k+1} \alpha_{k+1} \rho \la \Rb\xi(y_{k+1}), u - u_{k+1} \ra,
	\label{eq:Lm2}
	\end{align}	
	where $\E_{k+1}$ denotes the expectation conditioned on all the randomness up to step $k$. 
\end{lemma}

\begin{proof}
	First, for any $u\in Q$, by assumption~\ref{A:RIO}, 
	\begin{align}
	\mE_{k+1} \alpha_{k+1}\la \tnf(y_{k+1}), u_k - u\ra & \stackrel{~\eqref{eq:IODef}}{=} \mE_{k+1}\alpha_{k+1} \rho \la \Rb(\Rf^T\nabla f(y_{k+1}) + \xi(y_{k+1})), u_k - u \ra \notag \\
	& \stackrel{~\eqref{eq:IOProp1}}{=} \alpha_{k+1} \la \nabla f(y_{k+1}), u_k - u \ra \notag\\
	&\hspace{1.4em}+ \mE_{k+1}\alpha_{k+1} \rho \la \Rb\xi(y_{k+1}), u_k - u \ra.   
	\label{eq:Lm2Pr1}
	\end{align}	
	Taking conditional expectation $\mE_{k+1}$ in ~\eqref{eq:Lm1} of Lemma~\ref{Lm:1} and using ~\eqref{eq:Lm2Pr1}, we obtain the statement of the Lemma.	
\end{proof}

\begin{lemma}
	Let the sequences $\{x_k, y_k, u_k, \alpha_k, A_k \}$, $k\geq 0$ be generated by \ag{A}lgorithm~\ref{Alg:UARM}. Then,  for all $u \in Q$, it holds that
	\begin{align}
	A_{k+1} \E_{k+1}f(x_{k+1}) - A_{k}f(x_{k}) \leq & \,\,\alpha_{k+1} \left( f(y_{k+1}) + \la \nabla f(y_{k+1}), u - y_{k+1} \ra \right) + V[u_k](u)  \notag \\ 
	&- \E_{k+1}V[u_{k+1}](u) \notag\\
	&+ \E_{k+1} \alpha_{k+1} \rho \la \Rb\xi(y_{k+1}), u - u_{k+1} \ra.
	\label{eq:Lm3}
	\end{align}
	\label{Lm:3}
\end{lemma}

\begin{proof}
	For any $u \in Q$,
	\begin{align}
	\alpha_{k+1}\la \nabla f(y_{k+1}), y_{k+1} - u\ra &= \alpha_{k+1}\la \nabla f(y_{k+1}), y_{k+1} - u_k \ra + \alpha_{k+1}\la \nabla f(y_{k+1}), u_k - u\ra \notag \\
	&\hspace{-5.5em}\stackrel{~\eqref{eq:alpQuadEq},~\eqref{eq:ykp1Def}}{=} A_{k}\la \nabla f(y_{k+1}), x_k - y_{k+1} \ra + \alpha_{k+1}\la \nabla f(y_{k+1}), u_k - u\ra \notag \\
	&\hspace{-5.5em}\stackrel{\text{conv-ty}}{\leq} A_{k}\left(f(x_{k}) - f(y_{k+1}) \right) + \alpha_{k+1}\la \nabla f(y_{k+1}), u_k - u\ra \notag \\
	&\hspace{-5em}\stackrel{~\eqref{eq:Lm2}}{\leq} A_{k}\left(f(x_{k}) - f(y_{k+1}) \right) + A_{k+1}(f(y_{k+1}) - \E_{k+1}f(x_{k+1})) + \notag \\
	&\hspace{-3.5em} + V[u_k](u) - \E_{k+1}V[u_{k+1}](u) + 	\E_{k+1} \alpha_{k+1} \rho \la \Rb\xi(y_{k+1}), u - u_{k+1} \ra\notag \\
	&\hspace{-5em}  \stackrel{~\eqref{eq:alpQuadEq}}{=} \alpha_{k+1} f(y_{k+1}) + A_{k}f(x_{k}) - A_{k+1} \E_{k+1}f(x_{k+1}) + V[u_k](u) \notag \\
	&  \hspace{-3.5em}- \E_{k+1}V[u_{k+1}](u)  + 	\E_{k+1} \alpha_{k+1} \rho \la \Rb\xi(y_{k+1}), u - u_{k+1} \ra.
	\end{align}
	Rearranging terms, we obtain the statement of the Lemma.
\end{proof}

\begin{theorem}
	Let the assumptions~\ref{A:RIO}, \ref{A:RPM}, \ref{A:S} hold. Let the sequences $\{x_k, y_k, u_k, \alpha_k, A_k \}$, $k\geq 0$ be generated by \ag{A}lgorithm~\ref{Alg:UARM}. Let $f_*$ be the optimal objective value and $x_*$ be an optimal point in problem~\eqref{eq:PrSt}. Denote 
	\begin{equation}
	P_0^2 = A_0(f(x_0)-f_*) + V[u_0](x_*).
	\label{eq:P0Def}
	\end{equation}
	\begin{enumerate}
		\item If the oracle error $\xi(x)$ in ~\eqref{eq:IODef} can be controlled and, on each iteration, the error level $\delta$ in ~\eqref{eq:IOProp2} satisfies
		\begin{equation}
		\delta \leq \frac{P_0}{4\rho A_k},
		\label{eq:deltaContr}
		\end{equation}
		then, for all $k \geq 1$,
		$$
		\mE f(x_{k}) - f_*  \leq \frac{3P_0^2}{2A_{k}},
		$$
		where $\mE$ denotes the expectation with respect to all the randomness \ag{(}up to step $k$ \ag{in this case)}.
		\item If the oracle error $\xi(x)$ in ~\eqref{eq:IODef} can not be controlled, then, for all $k \geq 1$,
		$$
		\mE f(x_{k}) - f_*  \leq \frac{2P_0^2}{A_{k}} + 4 A_k \rho^2 \delta^2. 
		$$
	\end{enumerate}
	\label{Th:1}
\end{theorem}

\begin{proof}
	Let us change the counter in \ag{Lemma}~\ref{Lm:3} from $k$ to $i$, fix $u=x_*$, take the full expectation in each inequality for $i=0,...,k-1$ and sum all the inequalities for $i=0,...,k-1$. Then,  
	\begin{align}
	A_{k} \mE f(x_{k}) - A_{0}f(x_{0}) \leq & \,\, \sum_{i=0}^{k-1}\alpha_{i+1} \mE \left( f(y_{i+1}) + \la \nabla f(y_{i+1}), x_* - y_{i+1} \ra \right) + V[u_0](x_*) \notag \\ 
	&- \mE V[u_{k}](x_*) + \sum_{i=0}^{k-1} \mE \alpha_{i+1} \rho \la \Rb\xi(y_{i+1}), x_* - u_{i+1} \ra \notag \\
	& \hspace{-6em} \stackrel{\text{conv-ty},~\eqref{eq:alpQuadEq},~\eqref{eq:IOProp2}}{\leq} (A_k-A_0)f(x_*) + V[u_0](x_*)\notag\\
	&- \mE V[u_{k}](x_*) + \sum_{i=0}^{k-1} \alpha_{i+1} \rho \delta \mE  \|x_* - u_{i+1}\|_E \notag .
	\end{align}
	Rearranging terms and using ~\eqref{eq:P0Def}, we obtain, for all $k \geq 1$,
	\begin{align}
	0 \leq A_{k} \left(\mE f(x_{k}) - f_* \right) & \leq 
	P_0^2- \mE V[u_{k}](x_*) + \rho \delta \sum_{i=0}^{k-1} \alpha_{i+1} \mE R_{i+1},
	\label{eq:TelSumRes}
	\end{align}
	where we denoted $R_i = \|u_{i}-x_*\|_E$, $i\geq 0$. 
	
	1. We first prove the first statement of the Theorem. 
	We have
	\begin{equation}
	\frac12 R_0^2 = \frac12 \|x_*-u_{0}\|_E^2 \stackrel{~\eqref{eq:BFLowBound}}{\leq} V[u_0](x_*) \stackrel{~\eqref{eq:P0Def}}{\leq} P_0^2.
	\label{eq:}
	\end{equation}
	Hence, $\mE R_0 = R_0 \leq P_0 \sqrt{2} \leq 2 P_0$. Let $\mE R_i \leq 2 P_0$, for all $i=0,...,k-1$. Let us prove that $\mE R_{k} \leq 2 P_0$.
	By convexity of square function, we obtain
	\begin{align}
	\frac12 \left(\mE R_{k}\right)^2 \leq \frac12 \mE R_{k}^2 \stackrel{~\eqref{eq:BFLowBound}}{\leq} \mE V[u_{k}](x_*) & \stackrel{~\eqref{eq:TelSumRes}}{\leq} P_0^2 + \rho \delta \sum_{i=0}^{k-2} \alpha_{i+1} 2P_0 + \alpha_{k} \rho \delta \mE R_{k} \notag \\
	& \stackrel{~\eqref{eq:alpQuadEq}}{=} P_0^2 + 2 \rho \delta P_0 (A_{k-1}-A_0) + \alpha_{k} \rho \delta \mE R_{k} \notag \\
	& \leq P_0^2 + 2 \rho \delta P_0 A_{k} + \alpha_{k} \rho \delta \mE R_{k} .
	\label{eq:ERkp1ExpIneq}
	\end{align}
	Since $\alpha_k \leq A_k$, $k \geq 0$, by the choice of $\delta$ ~\eqref{eq:deltaContr}, we have 
	$2 \rho \delta P_0 A_k  \leq \frac{P_0^2}{2}$ and $\alpha_{k} \rho \delta \leq A_{k} \rho \delta \leq  \frac{P_0}{4}$.
	So, we obtain an inequality for $\mE R_{k}$
	$$
	\frac12 \left(\mE R_{k}\right)^2 \leq \frac{3P_0^2}{2} +  \frac{P_0}{4} \mE R_{k}.
	$$
	Solving this quadratic inequality in $\mE R_{k}$, we obtain
	$$
	\mE R_{k} \leq \frac{P_0}{4} + \sqrt{\frac{P_0^2}{16}+3P_0^2} = 2P_0.
	$$
	Thus, by induction, we have that, for all $k \geq 0$, $\mE R_k \leq 2 P_0$. Using the bounds $\mE R_i \leq 2 P_0$, for all $i=0,...,k$, we obtain
	\begin{align*}
	A_{k} \left(\mE f(x_{k}) - f_* \right) &\stackrel{~\eqref{eq:TelSumRes}}{\leq} P_0^2 + \rho \delta \sum_{i=0}^{k-1} \alpha_{i+1}\mE R_i\\ & \stackrel{~\eqref{eq:alpQuadEq},~\eqref{eq:deltaContr}}{\leq}  P_0^2 + \rho \frac{P_0}{4\rho A_k} \cdot (A_k-A_0) \cdot  2 P_0 \leq \frac{3P_0^2}{2}.
	\end{align*}
	This finishes the proof of the first statement of the Theorem.
	
	2. Now we prove the second statement of the Theorem.
	First, from \eqref{eq:TelSumRes} for $k=1$, we have
	\begin{align}
	\frac12 \left(\mE R_{1}\right)^2 \leq \frac12 \mE R_{1}^2 \stackrel{\eqref{eq:BFLowBound}}{\leq} \mE V[u_{1}](x_*) & \stackrel{\eqref{eq:TelSumRes}}{\leq} P_0^2 + \rho \delta \alpha_{1} \mE R_{1}.
	\notag
	\end{align} 
	Solving this inequality in $\mE R_{1}$, we obtain 
	\begin{equation}
	\mE R_{1} \leq  \rho \delta \alpha_{1}+ \sqrt{(\rho \delta\alpha_{1} )^2+2P_0^2} \leq 2  \rho \delta \alpha_{1} +  P_0  \sqrt{2},
	\end{equation}
	where we used that, for any $a,b \geq 0$, $\sqrt{a^2+b^2}\leq a +b$.
	Then,
	\begin{align}
	P_0^2 + \rho \delta \alpha_{1} \mE R_{1}  & \leq  P_0^2  + 2 (\rho \delta \alpha_1)^2 +  \rho \delta \alpha_1 P_0\sqrt{2}  \leq \left(P_0 + \rho \delta \sqrt{2} (A_{1}-A_0) \right)^2. \notag 
	\end{align}
	Thus, we have proved that the inequality
	\begin{equation}
	P_0^2 + \rho \delta \sum_{i=0}^{k-2}\alpha_{i+1} \mE R_{i+1} \leq \left(P_0 + \rho \delta \sqrt{2} (A_{k-1}-A_0)\right)^2
	\label{eq:Th1Pr1}
	\end{equation}
	holds for $k=2$. Let us assume that it holds for some $k$ and prove that it holds for $k+1$. We have
	\begin{align}
	\frac12 \left(\mE R_{k}\right)^2 \leq \frac12 \mE R_{k}^2 \stackrel{\eqref{eq:BFLowBound}}{\leq} \mE V[u_{k}](x_*) & \stackrel{\eqref{eq:TelSumRes}}{\leq} P_0^2 + \rho \delta \sum_{i=0}^{k-2} \alpha_{i+1} \mE R_{i+1} + \alpha_{k} \rho \delta \mE R_{k} \notag \\
	& \stackrel{\eqref{eq:Th1Pr1}}{\leq} \left(P_0 + \rho \delta \sqrt{2} (A_{k-1}-A_0)\right)^2 + \alpha_{k} \rho \delta \mE R_{k}. \notag 
	\end{align}
	Solving this quadratic inequality in $\mE R_{k}$, we obtain
	\begin{align}
	\mE R_{k} & \leq \alpha_{k} \rho \delta + \sqrt{(\alpha_{k} \rho \delta)^2+2\left(P_0 + \rho \delta \sqrt{2} (A_{k-1}-A_0)\right)^2} \notag \\
	& \leq 2 \alpha_{k} \rho \delta +  \left(P_0 + \rho \delta \sqrt{2} (A_{k-1}-A_0)\right) \sqrt{2},
	\label{eq:Th1Pr2}
	\end{align}
	where we used that, for any $a,b \geq 0$, $\sqrt{a^2+b^2}\leq a +b$. Further,
	\begin{align}
	P_0^2 + \rho \delta \sum_{i=0}^{k-1}\alpha_{i+1} \mE R_{i+1}  & \stackrel{\eqref{eq:Th1Pr1}}{\leq} \left(P_0 + \rho \delta \sqrt{2} (A_{k-1}-A_0)\right)^2 + \rho \delta \alpha_k \mE R_{k} \notag \\
	& \hspace{-9em} \stackrel{\eqref{eq:Th1Pr2}}{\leq} \left(P_0 + \rho \delta \sqrt{2} (A_{k-1}-A_0)\right)^2  + 2 (\rho \delta \alpha_k)^2 +  \rho \delta \alpha_k \left(P_0 + \rho \delta \sqrt{2} (A_{k-1}-A_0)\right) \sqrt{2} \notag \\
	& \hspace{-8.4em} \leq \left(P_0 + \rho \delta \sqrt{2} (A_{k-1}-A_0) + \rho \delta \alpha_k \sqrt{2} \right)^2 = \left(P_0 + \rho \delta \sqrt{2} (A_{k}-A_0) \right)^2, \notag 
	\end{align}
	which is \eqref{eq:Th1Pr1} for $k+1$.
	Using this inequality, we obtain 
	\begin{align*}
	A_{k} \left(\mE f(x_{k}) - f_* \right) & \stackrel{\eqref{eq:TelSumRes}}{\leq} P_0^2 + \rho \delta \sum_{i=0}^{k-1} \alpha_{i+1} \mE R_{i+1} \leq \left(P_0 + \rho \delta \sqrt{2} (A_{k}-A_0) \right)^2 \\
	&\hspace{0.6em}\leq 2P_0^2 + 4 \rho^2 \delta^2 A_k^2,
	\end{align*}
	which finishes the proof of the Theorem.

\end{proof}

Let us now estimate the growth rate of the sequence $A_k$, $k \geq 0$, which will give the rate of convergence for \ref{Alg:UARM}.
\begin{lemma}
	\label{Lm:AkGrowth}
	Let the sequence $\{A_k \}$, $k\geq 0$ be generated by \ag{Algorithm} \ref{Alg:UARM}. Then,  for all $k \geq 1$ it holds that
	\begin{equation}
	\frac{(k-1+2\rho)^2}{4\rho^2} \leq A_k \leq \frac{(k-1+2\rho)^2}{\rho^2}.
	\label{eq:AkGrowth}
	\end{equation}	
\end{lemma}

\begin{proof}
	As we showed in \ref{Lm:ConvComb}, $\alpha_1=\frac{1}{\rho}$ and, hence, $A_1=\alpha_0+\alpha_1=1$. Thus, \eqref{eq:AkGrowth} holds for $k=1$. Let us assume that \eqref{eq:AkGrowth} holds for some $k\geq 1$ and prove that it holds also for $k+1$.
	From \eqref{eq:alpQuadEq},  we have a quadratic equation for $\alpha_{k+1}$
	\begin{equation*}
	\rho^2\alpha_{k+1}^2 - \alpha_{k+1} - A_k = 0.
	\end{equation*}
	Since we need to take the largest root, we obtain,
	\begin{align}
	\alpha_{k+1} & = \frac{1 + \sqrt{\uprule 1 + 4\rho^2A_{k}}}{2\rho^2} = \frac{1}{2\rho^2} + \sqrt{\frac{1}{4\rho^4} + \frac{A_{k}}{\rho^2}} \geq 
	\frac{1}{2\rho^2} + \sqrt{\frac{A_{k}}{\rho^2}} \notag \\
	&\geq
	\frac{1}{2\rho^2} + \frac{k-1+2\rho}{2\rho^2} =
	\frac{k+2\rho}{2\rho^2}, \notag
	\end{align}
	where we used the induction assumption that \eqref{eq:AkGrowth} holds for $k$.
	On the other hand, 
	\begin{align}
	\alpha_{k+1} & = \frac{1}{2\rho^2} + \sqrt{\frac{1}{4\rho^4} + \frac{A_{k}}{\rho^2}} \leq 
	\frac{1}{\rho^2} + \sqrt{\frac{A_{k}}{\rho^2}} \notag \\
	&\leq
	\frac{1}{\rho^2} + \frac{k-1+2\rho}{\rho^2} =
	\frac{k+2\rho}{\rho^2}, \notag
	\end{align}
	where we used inequality $\sqrt{a+b} \leq \sqrt{a}+\sqrt{b}$, $a,b \geq 0$.
	Using the obtained inequalities for $\alpha_{k+1}$, from \eqref{eq:alpQuadEq} and \eqref{eq:AkGrowth} for $k$, we get
	\begin{equation*}
	A_{k+1} = A_k + \alpha_{k+1} \geq \frac{(k-1+2\rho)^2}{4\rho^2} + \frac{k+2\rho}{2\rho^2} \geq \frac{(k+2\rho)^2}{4\rho^2}
	\end{equation*}
	and
	\begin{equation*}
	A_{k+1} = A_k + \alpha_{k+1} \leq \frac{(k-1+2\rho)^2}{\rho^2} + \frac{k+2\rho}{\rho^2} \leq \frac{(k+2\rho)^2}{\rho^2}.
	\end{equation*}
	In the last inequality we used that $k \geq 1$, $\rho \geq 0$.
\end{proof}

\begin{remark}
	\label{Rm:Rate}
	According to \ag{Theorem}~\ref{Th:1}, if the desired accuracy of the solution is $\e$, i.e. the goal is to find such $\hat{x} \in Q$ that $\mE f(\hat{x}) - f_* \leq \e$, then the \ag{Algorithm}~\ref{Alg:UARM} should be stopped when $\frac{3P_0^2}{2A_k} \leq \e$. Then $\frac{1}{A_k} \leq \frac{2\e}{3 P_0^2}$ and the oracle error level $\delta$ should satisfy $\delta \leq \frac{P_0}{4 \rho A_k} \leq \frac{\e}{6 \rho P_0}$.
	
	From \ref{Lm:AkGrowth}, we obtain that $\frac{3P_0^2}{2A_k} \leq \e$ holds when $k$ is the smallest integer satisfying 
	$$
	\frac{(k-1+2\rho)^2}{4\rho^2} \geq \frac{3P_0^2}{2 \e}.
	$$
	This means that, to obtain an $\e$-solution, it is enough to choose
	\begin{equation}
	k = \max\left\{\left\lceil \rho \sqrt{\frac{6P_0^2}{\e}} + 1 -2 \rho \right\rceil,0\right\}.\label{eq:k_steps}
	\end{equation}
	Note that this dependence on $\e$ means that the proposed method is accelerated.
\end{remark}

\section{Extension for strongly convex functions}

In this section, we consider strongly-convex objective functions and show how the restart technique can be used to obtain faster rates of convergence under this additional assumption. 

\begin{assumption}[Strong convexity]\label{A:Strong_convexity}
	Assume that the function $f(x)$ is strongly convex:
	\begin{equation}
	f(x) \geq f(y) + \left\langle\nabla f(y), x-y\right\rangle + \frac{\mu}{2}\norm{x - y}_E^2, \quad \forall x, y \in Q,
	\end{equation}
\end{assumption}
where the constant $\mu > 0$. We also introduce additional assumption on the Bregman divergence $V[y](x)$:
\begin{assumption}[Prox-function bound]\label{A:Bregman_divergence}
	We assume that the function $d(x)$ satisfies conditions $0 = \argmin_{x\in Q} d(x)$ and $d(0) = 0$. Also for a fixed point $x_*$ and random point $y \in Q^0$, if $\mE \norm{y - x_*}^2_E \leq R^2$ then
	\begin{equation}
	\mE \, d\left(\frac{y - x_*}{R}\right) \leq \frac{V^2}{2}.
	\end{equation}
	
\end{assumption}
Some examples of prox-functions satisfying this assumption can be found in \cite{juditsky2014deterministic}, the simplest being the squared Euclidean norm. The following algorithm is obtained by restarting \ag{Algorithm}~\ref{Alg:UARM}. As in \ag{Algorithm}~\ref{Th:1} our algorithm has two variations. For simplicity we introduce a constant $\hat{C}$ which is equal to $0$ if the oracle error $\xi(x)$ in \eqref{eq:IODef} can be controlled and $\hat{C} = 4(k+2\rho)^2 \delta^2$ if the oracle error $\xi(x)$ in \eqref{eq:IODef} can not be controlled. Variable $k$ is a constant which is predefined in \ag{\eqref{eq:k_steps_A1}.} 

\begin{algorithm}[h!]
 	\caption{UARM for strongly convex functions}
	\label{Alg:UARM_strong}
 \SetKwInOut{Input}{Input}\SetKwInOut{Output}{Output}
 \KwIn{starting point $y_0 \in Q^0 = \otimes_{i=1}^n Q_i^0$, \at{number of iterations $N$}, prox-setup: $d(x)$, $V[u] (x)$, strong convexity constant $\mu > 0$, number $\e_{0}$ such that $f(y_0) - f_* \leq \e_{0}$}
 Set $j=0$. \\
 Set $A_0 = 1 - \frac{1}{\rho}$.\\
 Calculate \begin{equation}
			\label{eq:k_steps_A1}
			k = \max\left\{\left\lceil \rho \sqrt{16\left(1-\frac{1}{\rho} + \frac{V^2}{\mu}\right)} + 1 -2 \rho \right\rceil,0\right\}
			\end{equation}
 
 \If{the oracle error $\xi(x)$ in \eqref{eq:IODef} can be controlled}{Set $\hat{C} = 0$}
 
 \Else{Set $\hat{C} = 4(k+2\rho)^2 \delta^2$}

 \Repeat{\at{$j \leq N$}}{
  Set $j = j + 1$ \\
  Calculate $\e_{j} = \frac{\e_{j-1}}{2}$ \\
  \label{st:Alg1} Do $k$ steps of \ag{Algorithm}~\ref{Alg:UARM} with starting point $u_0 = y_{j-1}$, prox-setup:
			
			\begin{equation}
			\hat{d}(x) = \frac{2(\e_{j-1} + 2(1 - 2^{-j+1})\hat{C})}{\mu}d\left(\sqrt{\frac{\mu}{2(\e_{j-1}+ 2(1 - 2^{-j+1})\hat{C})}}
			(x - y_{j-1})\right),
			\label{auxiliary_prox}
			\end{equation}
		\at{corresponding Bregman divergence} $\hat{V}[u] (x)$, and output point $x_{j-1}$. \\
 Set $y_{j} = x_{j-1}$ \\
 }
 \KwOut{The point $y_{N}$.}
\end{algorithm}

\begin{theorem}
	Let the assumptions~\ref{A:RIO}, \ref{A:RPM}, \ref{A:S}, \ref{A:Strong_convexity},  \ref{A:Bregman_divergence} hold. Let the sequences $\{y_j, \e_j\}$, $j\geq 0$ be generated by \ag{A}lgorithm~\ref{Alg:UARM_strong}. Let $f_*$ be the optimal objective value and $x_*$ be an optimal point in problem~\eqref{eq:PrSt}.
	\begin{enumerate}
		\item For the second statement of \ag{Theorem}~\ref{Th:1} when the oracle error $\xi(x)$ in \eqref{eq:IODef} can not be controlled for all $j \geq 0$,
		\begin{align}
		\label{eq:Th2_2}
		\mE f(y_{j}) - f_* \leq \e_j + 8(1 - 2^{-j}) (k+2\rho)^2 \delta^2 = \e_0 2^{-j} + 8(1 - 2^{-j}) (k+2\rho)^2 \delta^2,
		\end{align}
		where $k$ is defined in \eqref{eq:k_steps_A1}.
		\item For the first statement of \ag{Theorem}~\ref{Th:1} when the oracle error $\xi(x)$ in \eqref{eq:IODef} can be controlled for all $j \geq 0$,
		\begin{align}
		\label{eq:Th2_1}
		\mE f(y_{j}) - f_* \leq \e_j = \e_0 2^{-j},
		\end{align}
		where $\mE$ denotes the expectation with respect to all the randomness up to step $j$. Also \ag{Algorithm} \ref{Alg:UARM_strong} finds an $\e$-solution of the function $f(x)$ in the number of \ref{Alg:UARM} steps not greater than
		\begin{equation*}
		\left\lceil \log_2\left(\frac{\e_0}{\e}\right) \right\rceil\max\left\{\left\lceil \rho \sqrt{16\left(1-\frac{1}{\rho} + \frac{V^2}{\mu}\right)} + 1 -2 \rho \right\rceil,0\right\}.
		\end{equation*}
	\end{enumerate}
	\label{Th:2}
\end{theorem}

\begin{proof}
	Let us denote $\hat{C} = 4(k+2\rho)^2 \delta^2$. We start our proof with a remark that if $d(x)$ is a prox-function then $a^2 d\left(\frac{1}{a}x\right)$ is also a prox-function, it is enough to show that $a ^2d\left(\frac{1}{a}x\right)$ is $1$-strong convex function. From definition
	\begin{equation}
	d(x) - d(y) - \langle\nabla d(y), x - y\rangle \geq \frac{1}{2}\norm{x - y}_E^2, \quad x \in Q, \quad y \in Q^0 := \otimes_{i=1}^n Q_i^0,
	\end{equation}
	obviously, we can get
	\begin{equation}
	a^2d\left(\frac{1}{a}x\right) - a^2d\left(\frac{1}{a}y\right) - \left\langle\nabla\left( a^2d\left(\frac{1}{a}y\right)\right), x - y\right\rangle \geq \frac{1}{2}\norm{x - y}_E^2,
	\end{equation}
	$x \in Q$, $y \in Q^0 := \otimes_{i=1}^n Q_i^0$.
	\at{Let us denote $\hat{P}_0^j =\sqrt{\left( A_0 + \frac{V^2}{\mu}\right)\e_{j}}$ for all $j=0,...,N$}, we can show following equality for $k$:
	\begin{align*}k &= \max\left\{\left\lceil \rho \sqrt{16\left(1-\frac{1}{\rho} + \frac{V^2}{\mu}\right)} + 1 -2 \rho \right\rceil,0\right\}\\
	&= \max\left\{\left\lceil \rho \sqrt{\frac{8\left(A_0 + \frac{V^2}{\mu}\right)\e_{j-1}}{\e_j}} + 1 -2 \rho \right\rceil,0\right\}\\
	&= \max\left\{\left\lceil \rho \sqrt{\frac{8(\hat{P}_0^{j-1})^2}{\e_j}} + 1 -2 \rho \right\rceil,0\right\}.\end{align*}
	Thus we have equivalent definition of $k$:
	\begin{align}
	k = \max\left\{\left\lceil \rho \sqrt{\frac{8(\hat{P}_0^{j-1})^2}{\e_j}} + 1 -2 \rho \right\rceil,0\right\}\label{eq:k_steps_A2}
	.\end{align}
	1. We first prove the first statement of \ag{Theorem}~\ref{Th:2}.
	For $j = 0$ inequality \eqref{eq:Th2_2} follows from \ag{Algorithm}~\ref{Alg:UARM_strong} input. Let us assume that this is true for $j$ and prove that it holds
	for $j+1$. Let us define $P_0^{j-1}$ equal to $P_0$ from the $j$'s call of \at{Algorithm} \ref{Alg:UARM}.
	From assumption~\ref{A:Strong_convexity}, induction, \at{and definition of $\hat{C}$}, we obtain 
	\begin{equation}
	\frac{\mu}{2}\mE\norm{y_{j} - x_*}_E^2 \leq \mE(f(y_{j}) - f_*) \leq \e_{j} + 2(1 - 2^{-j})\hat{C}.
	\label{eq:StrongConvOpt}
	\end{equation} 
	We can write following upper bound for $P_0^{j}$:
	\begin{align*}
	\mE(P_0^{j})^2 &= A_0\mE(f(y_j)-f_*) + \mE \hat{V}[y_j](x_*) \notag\\
	&=A_0\mE(f(y_j)-f_*) + \mE\left(\hat{d}(x_*) - \hat{d}(y_j) - \langle\nabla\hat{d}(y_j), x_* - y_j\rangle\right)
	\end{align*}
	From assumption~\ref{A:Bregman_divergence} and \eqref{auxiliary_prox}, we have $\hat{d}(y_j) = 0$ and $\nabla\hat{d}(y_j) = 0$. Therefore,
	\begin{align*}
	\mE(P_0^{j})^2&= A_0\mE(f(y_j)-f_*)\\ 
	&\hspace{1.5em}+ \frac{2(\e_j + 2(1 - 2^{-j})\hat{C})}{\mu}\mE d\left(\sqrt{\frac{\mu}{2(\e_{j}+ 2(1 - 2^{-j})\hat{C})}}(y_j - x_*)\right).
	\end{align*}
	In assumption~\ref{A:Bregman_divergence} we can take $R = \sqrt{\frac{2(\e_{j} + 2(1 - 2^{-j})\hat{C})}{\mu}}$ and get
	\begin{align*}
	\mE(P_0^{j})^2&\leq A_0\mE(f(y_j)-f_*) + \frac{(\e_{j}+ 2(1 - 2^{-j})\hat{C})V^2}{\mu} \notag\\
	&\leq\left(A_0 + \frac{V^2}{\mu}\right)(\e_{j}+ 2(1 - 2^{-j})\hat{C})\notag
	\end{align*}
	and hence $\mE(P_0^{j})^2 \leq (\hat{P}_0^{j})^2 + 2\left(A_0 + \frac{V^2}{\mu}\right)(1 - 2^{-j})\hat{C}$. Using Theorem \ref{Th:1} we obtain:
	
	\begin{align*}
	\mE f(y_{j+1}) - f_* 
	&\leq \frac{2\mE (P_0^{j})^2}{A_{k}} + 4 A_{k} \rho^2 \delta^2\\ 
	&\stackrel{\eqref{eq:AkGrowth}}{\leq} \frac{2\mE (P_0^{j})^2}{A_{k}} + \hat{C}\\ 
	&\leq \frac{2(\hat{P}_0^{j})^2}{A_{k}} + \frac{2\left(A_0 + \frac{V^2}{\mu}\right)(1 - 2^{-j})\hat{C}}{A_{k}} + \hat{C}.
	\end{align*}
	Using the second definition \eqref{eq:k_steps_A2} for $k$
	as well as in \ag{Remark}~\ref{Rm:Rate} with $P_0 = \hat{P}_0^{j}$ 
	we have:
	\begin{align*}
	\frac{2(\hat{P}_0^{j})^2}{A_{k}} \leq \e_{j+1}.
	\end{align*}
	Finally, we can conclude:
	\begin{align*}
	\mE f(y_{j+1}) - f_* &\leq \e_{j+1}+\frac{2\left(A_0 + \frac{V^2}{\mu}\right)(1 - 2^{-j})\hat{C}}{(\hat{P}_0^{j})^2}\e_{j+1}+ \hat{C}\\
	&=\e_{j+1}+\frac{2\left(A_0 + \frac{V^2}{\mu}\right)(1 - 2^{-j})\hat{C}}{\left( A_0 + \frac{V^2}{\mu}\right)\e_{j}}\e_{j+1}+ \hat{C}\\
	&=\e_{j+1}+2(1 - 2^{-j-1})\hat{C}.
	\end{align*}
	2. Obviously, the proof for the second statement of \ag{Theorem}~\ref{Th:2} is the same. The only difference is that $\hat{C} = 0$. It is enough to take $N = \left\lceil \log_2\left(\frac{\e_0}{\e}\right) \right\rceil$ in order to find an  $\e$-solution. Due to \eqref{eq:k_steps_A1} we can estimate total number of \ag{A}lgorithm~\ref{Alg:UARM} steps:
	\begin{align*}
	K_{T} &= \at{N  k} \\
	&= N\max\left\{\left\lceil \rho \sqrt{16\left(1-\frac{1}{\rho} + \frac{V^2}{\mu}\right)} + 1 -2 \rho \right\rceil,0\right\}\\
	&=\left\lceil \log_2\left(\frac{\e_0}{\e}\right) \right\rceil\max\left\{\left\lceil \rho \sqrt{16\left(1-\frac{1}{\rho} + \frac{V^2}{\mu}\right)} + 1 -2 \rho \right\rceil,0\right\}.\\
	\end{align*}
\end{proof}
Note that this dependence on $\e$ and $\mu$ means that the proposed method is accelerated.


\section{Examples of Applications}
\label{S:Appl}
In this section, we apply our general framework, which consists of assumptions~\ref{A:RIO}, \ref{A:RPM}, \ref{A:S}, UARM as listed in \ag{A}lgorithm~\ref{Alg:UARM} and convergence rate  in theorem~\ref{Th:1}, to obtain several particular algorithms and their convergence rate. We consider problem~\eqref{eq:PrSt} and, for each particular case, introduce a particular setup, which includes properties of the objective function $f$, available information about this function, properties of the feasible set $Q$. Based on each setup, we show how the Randomized Inexact Oracle is constructed and check that the assumptions~\ref{A:RIO}, \ref{A:RPM}, \ref{A:S} hold. 
Then, we obtain convergence rate guarantee for each particular algorithm as a corollary of theorem~\ref{Th:1}. Our examples include accelerated random directional search with inexact directional derivative, accelerated random block-coordinate descent with inexact block derivatives, accelerated random derivative-free directional search with inexact function values, accelerated random derivative-free block-coordinate descent with inexact function values. Accelerated random directional search and accelerated random derivative-free directional search were developed in \cite{nesterov2017random}, but for the case of exact directional derivatives and exact function values. Also, in the existing methods, a Gaussian random vector is used for randomization. Accelerated random block-coordinate descent was introduced in \cite{nesterov2012efficiency} and further developed in by several authors (see Introduction for the extended review). Existing methods of this type use exact information on the block derivatives and also only Euclidean proximal setup. In the contrast, our algorithm works with inexact derivatives and is able to work with entropy proximal setup. To the best of our knowledge, our accelerated random derivative-free block-coordinate descent with inexact function values is new. This method also can work with entropy proximal setup.

\subsection{Accelerated Random Directional Search}
\label{S:RDS}
In this subsection, we introduce accelerated random directional search with inexact directional derivative for unconstrained problems with Euclidean proximal setup.
We assume that, for all $i=1,...,n$, $Q_i=E_i=\R$,  $\|x^{(i)}\|_i^2=(x^{(i)})^2$, $x^{(i)} \in E_i$, $d_i(x^{(i)})=\frac12(x^{(i)})^2$, $x^{(i)} \in E_i$ and, hence, $V_i[z^{(i)}](x^{(i)}) = \frac12(x^{(i)}-z^{(i)})^2$, $x^{(i)}, z^{(i)} \in E_i$. Thus, $Q=E=\R^n$. Further, we assume that $f$ in \eqref{eq:PrSt} has $L$-Lipschitz-continuous gradient with respect to Euclidean norm, i.e. 
\begin{equation}
f(x) \leq f(y) + \la \nabla f(y), x-y \ra + \frac{L}{2}\|x-y\|_2^2, \quad x,y \in E.
\label{eq:fLipSm}
\end{equation}
We set $\beta_i=L$, $i=1,...,n$. Then, by definitions in subsection~\ref{S:Not}, we have $\|x\|_E^2 = L \|x\|_2^2$, $x \in E$, $d(x) = \frac{L}{2}\|x\|_2^2=\frac{1}{2}\|x\|_E^2$, $x \in E$, $V[z](x)=\frac{L}{2}\|x-z\|_2^2=\frac{1}{2}\|x-z\|_E^2$, $x, z \in E$. Also, we have $\|g\|_{E,*}^2 = L^{-1}\|g\|_{2}^2$, $g \in E^*$.

We assume that, at any point $x \in E$, one can calculate an inexact derivative of $f$ in a direction $e \in E$
$$
\tilde{f}'(x,e) = \la \nabla f(x), e \ra + \xi(x), 
$$ 
where $e$ is a random vector uniformly distributed on the Euclidean sphere of radius 1, i.e. $\Sp_2(1):=\{s \in \R^n: \|s\|_2=1\}$, and the directional derivative error $\xi(x) \in \R$ is uniformly bounded in absolute value by error level $\Delta$, i.e. $|\xi(x)| \leq \Delta$, $x \in E$. Since we are in the Euclidean setting, we consider $e$ also as an element of $E^*$. We use $n (\la \nabla f(x), e \ra + \xi(x))e$ as Randomized Inexact Oracle.

Let us check the assumptions stated in subsection~\ref{S:PrSt&Asmpt}.

\textbf{Randomized Inexact Oracle.}
In this setting, we have $\rho = n$, $H = \R$, $\Rf^T: E^* \to \R$ is given by $\Rf^T g = \la g, e \ra$, $g \in E^*$, $\Rb: \R \to E^*$ is given by $\Rb t = te$, $t \in \R$. Thus,
$$
\tnf(x) = n (\la \nabla f(x), e \ra + \xi(x))e.
$$
One can prove that $\E_e n \la \nabla f(x), e \ra e =  n\E_eee^T\nabla f(x) = \nabla f(x)$, $x \in E$, and, thus, \eqref{eq:IOProp1} holds. Also, for all $x \in E$, we have $\|\Rb \xi(x)\|_{E,*} = \frac{1}{\sqrt{L}}\|\xi(x)e\|_2 \leq \frac{\Delta}{\sqrt{L}}$, which proves \eqref{eq:IOProp2} if we take $\delta = \frac{\Delta}{\sqrt{L}}$.

\textbf{Regularity of Prox-Mapping.} Substituting particular choice of $Q$, $V[u](x)$, $\tnf(x)$ in \eqref{eq:uPlusDef}, we obtain
\begin{align*}
u_+ &= \arg \min_{x \in \R^n} \left\{\frac{L}{2}\|x-u\|_2^2+ \alpha \la n (\la \nabla f(y), e \ra + \xi(y))e,x \ra\right\}\\
& = u - \frac{\alpha n}{L} (\la \nabla f(y), e \ra + \xi(y))e.
\end{align*}
Hence, since $\la e, e \ra = 1$, we have
\begin{align}
\la \Rb\Rf^T\nabla f(y), u - u_+ \ra & = \left\la \la \nabla f(y), e \ra e, \frac{\alpha n}{L} (\la \nabla f(y), e \ra + \xi(y))e\right\ra  \notag \\
& =  \la \nabla f(y), e \ra \la e , e\ra \frac{\alpha n}{L} (\la \nabla f(y), e \ra + \xi(y)) \notag \\
& = \left\la \nabla f(y),  \frac{\alpha n}{L} (\la \nabla f(y), e \ra + \xi(y))e\right\ra \notag \\
& = \la \nabla f(y), u - u_+ \ra,  \notag
\end{align}
which proves \eqref{eq:RPMA}.

\textbf{Smoothness.}
By definition of $\|\cdot\|_E$ and \eqref{eq:fLipSm}, we have 
$$
f(x) \leq f(y) + \la \nabla f(y), x - y \ra + \frac{L}{2}\|x-y\|_2^2 = f(y) + \la \nabla f(y), x - y \ra + \frac{1}{2}\|x-y\|_E^2, \quad x, y \in E
$$
and \eqref{eq:NfLipAsm} holds.

We have checked that all the assumptions listed in subsection~\ref{S:PrSt&Asmpt} hold. Thus, we can obtain the following convergence rate result for random directional search as a corollary of \ag{Theorem}~\ref{Th:1} and \ag{Lemma}~\ref{Lm:AkGrowth}. 
\begin{corollary}
	\label{Cor:UARMDirSearch}
	Let Algorithm~\ref{Alg:UARM} with $\tnf(x) = n (\la \nabla f(x), e \ra + \xi(x))e$, where $e$ is random and uniformly distributed over the Euclidean sphere of radius 1, be applied to Problem \eqref{eq:PrSt} in the setting of this subsection. Let $f_*$ be the optimal objective value and $x_*$ be an optimal point in Problem \eqref{eq:PrSt}. Assume that directional derivative error $\xi(x)$ satisfies $|\xi(x)| \leq \Delta$, $x \in E$. Denote 
	$$
	P_0^2 = \left(1-\frac{1}{n}\right)(f(x_0)-f_*) + \frac{L}{2}\|u_0-x_*\|_2^2.
	$$
	\begin{enumerate}
		\item If the directional derivative error $\xi(x)$ can be controlled and, on each iteration, the error level $\Delta$ satisfies 
		$$
		\Delta \leq \frac{P_0\sqrt{L}}{4nA_k}, 
		$$
		then, for all $k \geq 1$,
		$$
		\mE f(x_{k}) - f_*  \leq \frac{6n^2P_0^2}{(k-1+2n)^2},
		$$
		where $\mE$ denotes the expectation with respect to all the randomness up to step $k$.
		\item If the directional derivative error $\xi(x)$ can not be controlled, then, for all $k \geq 1$,
		$$
		\mE f(x_{k}) - f_*  \leq \frac{8n^2P_0^2}{(k-1+2n)^2} + \frac{4}{L} (k-1+2n)^2 \Delta^2. 
		$$
	\end{enumerate}
\end{corollary}

\begin{remark}
	According to \ag{Remark}~\ref{Rm:Rate} and due to the relation $\delta = \frac{\Delta}{\sqrt{L}}$, we obtain that the error level $\Delta$ in the directional derivative should satisfy
	$$
	\Delta \leq \frac{\e \sqrt{L}}{6nP_0}. 
	$$
	At the same time, to obtain an $\e$-solution for Problem \eqref{eq:PrSt}, it is enough to choose
	$$
	k = \max\left\{\left\lceil n\sqrt{\frac{6P_0^2}{\e}} + 1 -2 n \right\rceil,0\right\}.
	$$
\end{remark}

\subsection{Accelerated Random Coordinate Descent}
\label{S:RCD}
In this subsection, we introduce accelerated random coordinate descent with inexact coordinate derivatives for problems with separable constraints and Euclidean proximal setup. We assume that, for all $i=1,...,n$, $E_i=\R$, $Q_i \subseteq E_i$ are closed and convex, $\|x^{(i)}\|_i^2=(x^{(i)})^2$, $x^{(i)} \in E_i$, $d_i(x^{(i)})=\frac12(x^{(i)})^2$, $x^{(i)} \in Q_i$, and, hence, $V_i[z^{(i)}](x^{(i)}) = \frac12(x^{(i)}-z^{(i)})^2$, $x^{(i)}, z^{(i)} \in Q_i$. Thus, $Q=\otimes_{i=1}^n Q_i$ has separable structure. 

Let us denote $e_i \in E$ the $i$-th coordinate vector. Then, for $i=1,...,n$, the $i$-th coordinate derivative of $f$ is $f_i'(x) = \la  \nabla f(x), e_i \ra $.
We assume that the gradient of $f$ in \eqref{eq:PrSt} is coordinate-wise Lipschitz continuous with constants $L_i$, $i=1,...,n$, i.e.
\begin{equation}
|f'_i(x+he_i) - f'_i(x)| \leq L_i |h|, \quad h \in \R, \quad i=1,...,n, \quad x \in Q.
\label{eq:fCoordLipSm}
\end{equation}
We set $\beta_i=L_i$, $i=1,...,n$. Then, by definitions in subsection~\ref{S:Not}, we have $\|x\|_E^2 = \sum_{i=1}^n L_i (x^{(i)})^2$, $x \in E$, $d(x) = \frac12\sum_{i=1}^n L_i(x^{(i)})^2$, $x \in Q$, $V[z](x)=\frac12\sum_{i=1}^n L_i (x^{(i)}-z^{(i)})^2$, $x, z \in Q$. Also, we have $\|g\|_{E,*}^2 = \sum_{i=1}^n L_i^{-1} (g^{(i)})^2$, $g \in E^*$. 

We assume that, at any point $x \in Q$, one can calculate an inexact coordinate derivative of $f$
$$
\tilde{f}_i'(x) = \la  \nabla f(x), e_i \ra + \xi(x), 
$$ 
where the coordinate $i$ is chosen from $i=1,...,n$ at random with uniform probability $\frac{1}{n}$, the coordinate derivative error $\xi(x) \in \R$ is uniformly bounded in absolute value by $\Delta$, i.e. $|\xi(x)| \leq \Delta$, $x \in Q$. Since we are in the Euclidean setting, we consider $e_i$ also as an element of $E^*$. We use $n (\la  \nabla f(x), e_i \ra + \xi(x))e_i$ as Randomized Inexact Oracle.

Let us check the assumptions stated in subsection~\ref{S:PrSt&Asmpt}.

\textbf{Randomized Inexact Oracle.}
In this setting, we have $\rho = n$, $H = E_i = \R$, $\Rf^T: E^* \to \R$ is given by $\Rf^T g = \la g, e_i \ra$, $g \in E^*$, $\Rb: \R \to E^*$ is given by $\Rb t = te_i$, $t \in \R$. Thus,
$$
\tnf(x) = n (\la \nabla f(x), e_i \ra + \xi(x))e_i, \quad x \in Q.
$$
One can prove that $\mE_i n \la \nabla f(x), e_i \ra e_i =  n\E_ie_ie_i^T\nabla f(x) = \nabla f(x)$, $x \in Q$, and, thus, \eqref{eq:IOProp1} holds. Also, for all $x \in Q$, we have $\|\Rb \xi(x)\|_{E,*} = \frac{1}{\sqrt{L_i}}|\xi(x)| \leq \frac{\Delta}{\sqrt{L_0}}$, where $L_0 = \min_{i=1,...,n}L_i$. This proves \eqref{eq:IOProp2} with $\delta = \frac{\Delta}{\sqrt{L_0}}$.

\textbf{Regularity of Prox-Mapping.} Separable structure of $Q$ and $V[u](x)$ means that the problem \eqref{eq:uPlusDef} boils down to $n$ independent problems of the form 
$$
u_+^{(j)} = \arg \min_{x^{(j)}\in Q_j} \left\{\frac{L_j}{2}(u^{(j)}-x^{(j)})^2+ \alpha \la \tnf(y), e_j \ra x^{(j)} \right\}, \quad j = 1,...,n.
$$
Since $\tnf(y)$ has only one, $i$-th, non-zero component, $\la \tnf(y), e_j \ra$ is zero for all $j \ne i$. Thus, $u-u_+$ has one, $i$-th, non-zero component and $\la  e_i, u - u_+ \ra e_i = u - u_+ $. 
Hence, 
\begin{align}
\la \Rb\Rf^T\nabla f(y), u - u_+ \ra & = \la \la \nabla f(y), e_i \ra e_i, u - u_+ \ra \notag \\
& = \la \nabla f(y), e_i \ra \la  e_i, u - u_+ \ra \notag \\
& = \la \nabla f(y), \la  e_i, u - u_+ \ra e_i \ra  \notag \\
& = \la \nabla f(y), u - u_+ \ra,  \notag
\end{align}
which proves \eqref{eq:RPMA}.

\textbf{Smoothness.}
By the standard reasoning, using \eqref{eq:fCoordLipSm}, one can prove that, for all $i=1,...,n$,
\begin{equation}
f(x+h e_i) \leq f(x) + h \la \nabla f(x),  e_i \ra  + \frac{L_i h^2}{2}, \quad h \in \R, \quad x \in Q.
\label{eq:fCoordLipSm2}
\end{equation}
Let $u,y \in Q$, $a \in \R$, and $x = y + a (u_+-u) \in Q$. As we have shown above, $u_+-u$ has only one, $i$-th, non-zero component. Hence, there exists $h \in \R$, such that $u_+-u = he_i $ and $x = y + a he_i $.
Thus, by definition of $\|\cdot\|_E$ and \eqref{eq:fCoordLipSm2}, we have 
\begin{align}
f(x)  = f(y + a he_i ) & \leq f(y) + a h \la \nabla f(y), e_i\ra    + \frac{L_i}{2}(a h)^2 \notag \\
& = f(y) + \la \nabla f(y), a he_i \ra + \frac{1}{2}\|a h e_i\|_E^2  \notag \\
& = f(y) + \la \nabla f(y), x - y \ra + \frac{1}{2}\|x-y\|_E^2. \notag 
\end{align}
This proves \eqref{eq:NfLipAsm}.

We have checked that all the assumptions listed in subsection~\ref{S:PrSt&Asmpt} hold. Thus, we can obtain the following convergence rate result for random coordinate descent as a corollary of \ref{Th:1} and \ref{Lm:AkGrowth}. 
\begin{corollary}
	\label{Cor:UARMCD}
	Let Algorithm~\ref{Alg:UARM} with $\tnf(x) = n (\la \nabla f(x), e_i \ra + \xi(x))e_i$, where $i$ is uniformly at random chosen from $1,...,n$, be applied to Problem \eqref{eq:PrSt} in the setting of this subsection. Let $f_*$ be the optimal objective value and $x_*$ be an optimal point in Problem \eqref{eq:PrSt}. Assume that coordinate derivative error $\xi(x)$ satisfies $|\xi(x)| \leq \Delta$, $x \in Q$. Denote 
	$$
	P_0^2 = \left(1-\frac{1}{n}\right)(f(x_0)-f_*) + \sum_{i=1}^n\frac{L_i}{2}(u_0^{(i)}-x_*^{(i)})^2.
	$$
	\begin{enumerate}
		\item If the coordinate derivative error $\xi(x)$ can be controlled and, on each iteration, the error level $\Delta$ satisfies 
		$$
		\Delta \leq \frac{P_0\sqrt{L_0}}{4nA_k}, 
		$$
		then, for all $k \geq 1$,
		$$
		\mE f(x_{k}) - f_*  \leq \frac{6n^2P_0^2}{(k-1+2n)^2},
		$$
		where $\mE$ denotes the expectation with respect to all the randomness up to step $k$.
		\item If the coordinate derivative error $\xi(x)$ can not be controlled, then, for all $k \geq 1$,
		$$
		\mE f(x_{k}) - f_*  \leq \frac{8n^2P_0^2}{(k-1+2n)^2} + \frac{4}{L_0} (k-1+2n)^2 \Delta^2. 
		$$
	\end{enumerate}
\end{corollary}

\begin{remark}
	According to \ag{Remark}~\ref{Rm:Rate} and due to the relation $\delta = \frac{\Delta}{\sqrt{L_0}}$, we obtain that the error level $\Delta$ in the coordinate derivative should satisfy
	$$
	\Delta \leq \frac{\e \sqrt{L_0}}{6nP_0} . 
	$$
	At the same time, to obtain an $\e$-solution for Problem \eqref{eq:PrSt}, it is enough to choose
	$$
	k = \max\left\{\left\lceil n\sqrt{\frac{6P_0^2}{\e}} + 1 -2 n \right\rceil,0\right\}.
	$$
\end{remark}

\subsection{Accelerated Random Block-Coordinate Descent}
\label{S:RBCD}
In this subsection, we consider two block-coordinate settings. The first one is the Euclidean, which is usually used in the literature for accelerated block-coordinate descent. The second one is the entropy, which, to the best of our knowledge, is analyzed in this context for the first time. We develop accelerated random block-coordinate descent with inexact block derivatives for problems with simple constraints in these two settings and their combination.

\textit{Euclidean setup.} We assume that, for all $i=1,...,n$, $E_i=\R^{p_i}$; $Q_i$ is a simple closed convex set; $\|x^{(i)}\|_i^2=\la B_i x^{(i)}, x^{(i)} \ra$, $x^{(i)} \in E_i$, where $B_i$ is symmetric positive semidefinite matrix; $d_i(x^{(i)})=\frac12\|x^{(i)}\|_i^2$, $x^{(i)} \in Q_i$, and, hence, $V_i[z^{(i)}](x^{(i)}) = \frac12\|x^{(i)}-z^{(i)}\|_i^2$, $x^{(i)}, z^{(i)} \in Q_i$.

\textit{Entropy setup.} We assume that, for all $i=1,...,n$, $E_i=\R^{p_i}$; $Q_i$ is standard simplex in $\R^{p_i}$, i.e., $Q_i = \{x^{(i)} \in \R^{p_i}_+: \sum_{j=1}^{p_i}[x^{(i)}]_j = 1\}$; $\|x^{(i)}\|_i=\|x^{(i)}\|_1=\sum_{j=1}^{p_i}|[x^{(i)}]_j|$, $x^{(i)} \in E_i$; $d_i(x^{(i)})=\sum_{j=1}^{p_i}[x^{(i)}]_j \ln [x^{(i)}]_j$, $x^{(i)} \in Q_i$, and, hence, $V_i[z^{(i)}](x^{(i)}) = \sum_{j=1}^{p_i}[x^{(i)}]_j \ln \frac{[x^{(i)}]_j}{[z^{(i)}]_j}$, $x^{(i)}, z^{(i)} \in Q_i$. 

Note that, in each block, one also can choose other proximal setups from \cite{ben-tal2015lectures}. Combination of different setups in different blocks is also possible, i.e., in one block it is possible to choose the Euclidean setup and in an another block one can choose the entropy setup.

Using operators $U_i$, $i=1,...,n$ defined in \eqref{eq:UDef}, for each $i=1,...,n$, the $i$-th block derivative of $f$ can be written as $f_i'(x) = U_i^T \nabla f(x)$.
We assume that the gradient of $f$ in \eqref{eq:PrSt} is block-wise Lipschitz continuous with constants $L_i$, $i=1,...,n$ with respect to chosen norms $\|\cdot \|_i$, i.e.
\begin{equation}
\|f'_i(x+U_ih^{(i)}) - f'_i(x)\|_{i,*} \leq L_i \|h^{(i)}\|_i, \quad h^{(i)} \in E_i, \quad i=1,...,n, \quad x \in Q.
\label{eq:fBlockLipSm}
\end{equation}
We set $\beta_i=L_i$, $i=1,...,n$. Then, by definitions in subsection~\ref{S:Not}, we have
\begin{equation*}
\|x\|_E^2 = \sum_{i=1}^n L_i \|x^{(i)}\|_i^2, \quad x \in E,
\end{equation*}
\begin{equation*}
d(x) = \sum_{i=1}^n L_id_i(x^{(i)}), \quad x \in Q,
\end{equation*}
\begin{equation*}
	V[z](x)=\sum_{i=1}^n L_i V_i[z^{(i)}](x^{(i)}), \quad x, z \in Q.
\end{equation*}
 Also, we have $\|g\|_{E,*}^2 = \sum_{i=1}^n L_i^{-1} \|g^{(i)}\|_{i,*}^2$, $g \in E^*$. 

We assume that, at any point $x \in Q$, one can calculate an inexact block derivative of $f$
$$
\tilde{f}_i'(x) = U_i^T \nabla f(x) + \xi(x),  
$$ 
where a block number $i$ is chosen from $1,...,n$ randomly uniformly, the block derivative error $\xi(x) \in E_i^*$ is uniformly bounded in norm by $\Delta$, i.e. $\|\xi(x)\|_{i,*} \leq \Delta$, $x \in Q$, $i=1,...,n$. As Randomized Inexact Oracle, we use $n \widetilde{U}_i(U_i^T \nabla f(x) + \xi(x))$, where $\widetilde{U}_i$ is defined in \eqref{eq:tUDef}.  

Let us check the assumptions stated in subsection~\ref{S:PrSt&Asmpt}.

\textbf{Randomized Inexact Oracle.}
In this setting, we have $\rho = n$, $H = E_i$, $\Rf^T: E^* \to E_i^*$ is given by $\Rf^T g = U_i^T g$, $g \in E^*$, $\Rb: E_i^* \to E^*$ is given by $\Rb g^{(i)} = \widetilde{U}_i g^{(i)}$, $g^{(i)} \in E_i^*$. Thus,
$$
\tnf(x) = n \widetilde{U}_i(U_i^T \nabla f(x) + \xi(x)), \quad x \in Q.
$$
Since $i \in R[1,n]$, one can prove that $\mE_i n \widetilde{U}_iU_i^T \nabla f(x) =  \nabla f(x)$, $x \in Q$, and, thus, \eqref{eq:IOProp1} holds. Also, for all $x \in Q$, we have $\|\Rb \xi(x)\|_{E,*} = \|\widetilde{U}_i \xi(x)\|_{E,*} = \frac{1}{\sqrt{L_i}}\|\xi(x)\|_{i,*} \leq \frac{\Delta}{\sqrt{L_0}}$, where $L_0 = \min_{i=1,...,n}L_i$. This proves \eqref{eq:IOProp2} with $\delta = \frac{\Delta}{\sqrt{L_0}}$. 

\textbf{Regularity of Prox-Mapping.} Separable structure of $Q$ and $V[u](x)$ means that the problem \eqref{eq:uPlusDef} boils down to $n$ independent problems of the form 
$$
u_+^{(j)} = \arg \min_{x^{(j)} \in Q_j} \left\{L_jV[u^{(j)}](x^{(j)})+ \alpha \la U_j^T\tnf(y),x^{(j)} \ra\right\}, \quad j = 1,...,n.
$$
Since $\tnf(y)$ has non-zero components only in the block $i$, $U_j^T\tnf(y)$ is zero for all $j \ne i$. Thus, $u-u_+$ has non-zero components only in the block $i$ and $U_i\widetilde{U}_i^T (u - u_+) = u - u_+$. 
Hence, 
\begin{align}
\la \Rb\Rf^T\nabla f(y), u - u_+ \ra & = \la \widetilde{U}_iU_i^T\nabla f(y), u - u_+ \ra \notag \\
& = \la \nabla f(y), U_i\widetilde{U}_i^T (u - u_+) \ra \notag \\
& = \la \nabla f(y), u - u_+ \ra,  \notag
\end{align}
which proves \eqref{eq:RPMA}.

\textbf{Smoothness.}
By the standard reasoning, using \eqref{eq:fBlockLipSm}, one can prove that, for all $i=1,...,n$,
\begin{equation}
f(x+U_ih^{(i)}) \leq f(x) + \la U_i^T\nabla f(x), h^{(i)} \ra + \frac{L_i}{2}\|h^{(i)}\|_i^2, \quad h^{(i)} \in E_i, \quad x \in Q.
\label{eq:fBlockLipSm2}
\end{equation}
Let $u,y \in Q$, $a \in \R$, and $x = y + a (u_+-u) \in Q$. As we have shown above, $u_+-u$ has non-zero components only in the block $i$. Hence, there exists $h^{(i)} \in E_i$, such that $u_+-u = U_i h^{(i)} $ and $x = y + a U_i h^{(i)} $.
Thus, by definition of $\|\cdot\|_E$ and \eqref{eq:fBlockLipSm2}, we have 
\begin{align}
f(x)  = f(y + a U_i h^{(i)}) & \leq f(y) + \la U_i^T\nabla f(y), a h^{(i)} \ra + \frac{L_i}{2}\|a h^{(i)}\|_i^2 \notag \\
& = f(y) + \la \nabla f(y), a U_i h^{(i)} \ra + \frac{1}{2}\|a U_i h^{(i)}\|_E^2  \notag \\
& = f(y) + \la \nabla f(y), x - y \ra + \frac{1}{2}\|x-y\|_E^2. \notag 
\end{align}
This proves \eqref{eq:NfLipAsm}.

We have checked that all the assumptions listed in subsection~\ref{S:PrSt&Asmpt} hold. Thus, we can obtain the following convergence rate result for random block-coordinate descent as a corollary of \ref{Th:1} and \ref{Lm:AkGrowth}. 
\begin{corollary}
	\label{Cor:UARMBlockCoord}
	Let \ref{Alg:UARM} with $\tnf(x) = n \widetilde{U}_i(U_i^T \nabla f(x) + \xi(x))$, where $i$ is uniformly at random chosen from $1,...,n$, be applied to Problem \eqref{eq:PrSt} in the setting of this subsection. Let $f_*$ be the optimal objective value and $x_*$ be an optimal point in Problem \eqref{eq:PrSt}. Assume that block derivative error $\xi(x)$ satisfies $|\xi(x)| \leq \Delta$, $x \in Q$. Denote 
	$$
	P_0^2 = \left(1-\frac{1}{n}\right)(f(x_0)-f_*) + V[u_0](x_*).
	$$
	\begin{enumerate}
		\item If the block derivative error $\xi(x)$ can be controlled and, on each iteration, the error level $\Delta$ satisfies 
		$$
		\Delta \leq \frac{P_0\sqrt{L_0}}{4nA_k}, 
		$$
		then,  for all $k \geq 1$,
		$$
		\mE f(x_{k}) - f_*  \leq \frac{6n^2P_0^2}{(k-1+2n)^2},
		$$
		where $\mE$ denotes the expectation with respect to all the randomness up to step $k$.
		\item If the block derivative error $\xi(x)$ can not be controlled, then, for all $k \geq 1$,
		$$
		\mE f(x_{k}) - f_*  \leq \frac{8n^2P_0^2}{(k-1+2n)^2} + \frac{4}{L_0} (k-1+2n)^2 \Delta^2. 
		$$
	\end{enumerate}
\end{corollary}

\begin{remark}
	According to \ag{Remark}~\ref{Rm:Rate} and due to the relation $\delta = \frac{\Delta}{\sqrt{L_0}}$, we obtain that the block derivative error $\Delta$ should satisfy
	$$
	\Delta \leq \frac{\e \sqrt{L_0}}{6nP_0}. 
	$$
	At the same time, to obtain an $\e$-solution for Problem \eqref{eq:PrSt}, it is enough to choose
	$$
	k = \max\left\{\left\lceil n\sqrt{\frac{6P_0^2}{\e}} + 1 -2 n \right\rceil,0\right\}.
	$$
\end{remark}

\subsection{Accelerated Random Derivative-Free Directional Search}
\label{S:RDFDS}
In this subsection, we consider the same setting as in subsection~\ref{S:RDS}, except for Randomized Inexact Oracle. Instead of directional derivative, we use here its finite-difference approximation.
We assume that, for all $i=1,...,n$, $Q_i=E_i=\R$,  $\|x^{(i)}\|_i=(x^{(i)})^2$, $x^{(i)} \in E_i$, $d_i(x^{(i)})=\frac12(x^{(i)})^2$, $x^{(i)} \in E_i$, and, hence, $V_i[z^{(i)}](x^{(i)}) = \frac12(x^{(i)}-z^{(i)})^2$, $x^{(i)}, z^{(i)} \in E_i$. Thus, $Q=E=\R^n$. Further, we assume that $f$ in \eqref{eq:PrSt} has $L$-Lipschitz-continuous gradient with respect to Euclidean norm, i.e. 
\begin{equation}
f(x) \leq f(y) + \la \nabla f(y), x-y \ra + \frac{L}{2}\|x-y\|_2^2, \quad x,y \in E.
\label{eq:fLipSmRDFDS}
\end{equation}
We set $\beta_i=L$, $i=1,...,n$. Then, by definitions in subsection~\ref{S:Not}, we have $\|x\|_E^2 = L \|x\|_2^2$, $x \in E$, $d(x) = \frac{L}{2}\|x\|_2^2=\frac{1}{2}\|x\|_E^2$, $x \in E$, $V[z](x)=\frac{L}{2}\|x-z\|_2^2=\frac{1}{2}\|x-z\|_E^2$, $x, z \in E$. Also, we have $\|g\|_{E,*}^2 = L^{-1}\|g\|_{2}^2$, $g \in E^*$.

We assume that, at any point $x \in E$, one can calculate an inexact value $\tf(x)$ of the function $f$, s.t. $|\tf(x)-f(x)|\leq \Delta$, $x \in E$. To approximate the gradient of $f$, we use 
$$
\tnf(x) = n\frac{\tf(x+\tau e)-\tf(x)}{\tau}e, 
$$ 
where $\tau > 0$ is small parameter, which will be chosen later, $e \in E$ is a random vector uniformly distributed on the Euclidean sphere of radius 1, i.e. on $\Sp_2(1):=\{s \in \R^n: \|s\|_2=1\}$.
Since, we are in the Euclidean setting, we consider $e$ also as an element of $E^*$. 

Let us check the assumptions stated in subsection~\ref{S:PrSt&Asmpt}.

\textbf{Randomized Inexact Oracle.}
First, let us show that the finite-difference approximation for the gradient of $f$ can be expressed in the form of \eqref{eq:IODef}. We have
\begin{align*}
\tnf(x) &= n\frac{\tf(x+\tau e)-\tf(x)}{\tau}e \\
& = n \left(\langle \nabla f(x), e\rangle  + \frac{1}{\tau}( \tf(x+\tau e)-\tf(x) - \tau \langle \nabla f(x), e\rangle)  \right)e.  \notag  
\end{align*}
Taking $\rho = n$, $H = \R$, $\Rf^T: E^* \to \R$ be given by $\Rf^T g = \la g, e \ra$, $g \in E^*$, $\Rb: \R \to E^*$ be given by $\Rb t = te$, $t \in \R$, we obtain
$$
\tnf(x) = n (\la \nabla f(x), e \ra + \xi(x))e,
$$
where $\xi(x) = \frac{1}{\tau}( \tf(x+\tau e)-\tf(x) - \tau \langle \nabla f(x), e\rangle)$.
One can prove that $\E_e n \la \nabla f(x), e \ra e =  n\E_eee^T\nabla f(x) = \nabla f(x)$, $x \in E$, and, thus, \eqref{eq:IOProp1} holds. It remains to prove \eqref{eq:IOProp2}, i.e., find $\delta$ s.t. for all $x \in E$, we have $\|\Rb \xi(x)\|_{E,*} \leq \delta$.
\begin{align}
\|\Rb \xi(x)\|_{E,*}  & =  \frac{1}{\sqrt{L}}\|\xi(x)e\|_2  = \frac{1}{\sqrt{L}}\left\|\frac{1}{\tau}( \tf(x+\tau e)-\tf(x) - \tau \langle \nabla f(x), e\rangle) e \right\|_2\notag  \\
& = \frac{1}{\sqrt{L}}\left\|\frac{1}{\tau}( \tf(x+\tau e) - f(x+\tau e) - (\tf(x) - f(x)) \right\|_2 \notag \\ 
& \hspace{1em} + \frac{1}{\sqrt{L}}\left\|(f(x+\tau e) - f(x) - \tau \langle \nabla f(x), e\rangle)) e \right\|_2 \notag \\
& \leq \frac{2\Delta}{\tau\sqrt{L}}+ \frac{\tau\sqrt{L}}{2}. \notag
\end{align}
Here we used that $|\tf(x)-f(x)|\leq \Delta$, $x \in E$ and \eqref{eq:fLipSmRDFDS}. So, we have that \eqref{eq:IOProp2} holds with $\delta = \frac{2\Delta}{\tau\sqrt{L}}+ \frac{\tau\sqrt{L}}{2}$. To balance both terms, we choose $\tau = 2 \sqrt{\frac{\Delta}{L}}$, which leads to equality $\delta = 2 \sqrt{\Delta}$.

\textbf{Regularity of Prox-Mapping.} This assumption can be checked in the same way as in subsection~\ref{S:RDS}.

\textbf{Smoothness.} This assumption can be checked in the same way as in subsection~\ref{S:RDS}.

We have checked that all the assumptions listed in subsection~\ref{S:PrSt&Asmpt} hold. Thus, we can obtain the following convergence rate result for random derivative-free directional search as a corollary of \ag{Theorem}~\ref{Th:1} and \ag{Lemma}~\ref{Lm:AkGrowth}. 
\begin{corollary}
	\label{Cor:UARMRDFDS}
	Let Algorithm~\ref{Alg:UARM} with $\tnf(x) = n\frac{\tf(x+\tau e)-\tf(x)}{\tau}e$, where $e$ is random and uniformly distributed over the Euclidean sphere of radius 1, be applied to Problem \eqref{eq:PrSt} in the setting of this subsection. Let $f_*$ be the optimal objective value and $x_*$ be an optimal point in Problem \eqref{eq:PrSt}. Assume that function value error $\tf(x)-f(x)$ satisfies $|\tf(x)-f(x)| \leq \Delta$, $x \in E$. Denote 
	$$
	P_0^2 = \left(1-\frac{1}{n}\right)(f(x_0)-f_*) + \frac{L}{2}\|u_0-x_*\|_2^2.
	$$
	\begin{enumerate}
		\item If the error in the value of the objective $f$ can be controlled and, on each iteration, the error level $\Delta$ satisfies 
		$$
		\Delta \leq \frac{P_0^2}{64n^2A_k^2}, 
		$$
		and $\tau = 2 \sqrt{\frac{\Delta}{L}}$ then, for all $k \geq 1$,
		$$
		\mE f(x_{k}) - f_*  \leq \frac{6n^2P_0^2}{(k-1+2n)^2},
		$$
		where $\mE$ denotes the expectation with respect to all the randomness up to step $k$.
		\item If the error in the value of the objective $f$ can not be controlled and $\tau = 2 \sqrt{\frac{\Delta}{L}}$, then, for all $k \geq 1$,
		$$
		\mE f(x_{k}) - f_*  \leq \frac{8n^2P_0^2}{(k-1+2n)^2} + 16 (k-1+2n)^2 L\Delta. 
		$$
	\end{enumerate}
\end{corollary}

\begin{remark}
	According to \ag{Remark}~\ref{Rm:Rate} and due to the relation $\delta = 2 \sqrt{\Delta}$, we obtain that the error level in the function value should satisfy
	$$
	\Delta \leq \frac{\e^2 }{144 n^2P_0^2}. 
	$$
	The parameter $\tau$ should satisfy
	$$
	\tau \leq \frac{\e }{6 n P_0 \sqrt{L}} .
	$$
	At the same time, to obtain an $\e$-solution for Problem \eqref{eq:PrSt}, it is enough to choose
	$$
	k = \max\left\{\left\lceil n\sqrt{\frac{6P_0^2}{\e}} + 1 -2 n \right\rceil,0\right\}.
	$$
\end{remark}

\subsection{Accelerated Random Derivative-Free Coordinate Descent}
\label{S:RDFCD}
In this subsection, we consider the same setting as in subsection~\ref{S:RCD}, except for Randomized Inexact Oracle. Instead of coordinate derivative, we use here its finite-difference approximation. We assume that, for all $i=1,...,n$, $E_i=\R$, $Q_i \subseteq E_i$ are closed and convex, $\|x^{(i)}\|_i=(x^{(i)})^2$, $x^{(i)} \in E_i$, $d_i(x^{(i)})=\frac12(x^{(i)})^2$, $x^{(i)} \in Q_i$, and, hence, $V_i[z^{(i)}](x^{(i)}) = \frac12(x^{(i)}-z^{(i)})^2$, $x^{(i)}, z^{(i)} \in Q_i$. Thus, $Q=\otimes_{i=1}^n Q_i$ has separable structure. 

Let us denote $e_i \in E$ the $i$-th coordinate vector. Then, for $i=1,...,n$, the $i$-th coordinate derivative of $f$ is $f_i'(x) = \la  \nabla f(x), e_i \ra $.
We assume that the gradient of $f$ in \eqref{eq:PrSt} is coordinate-wise Lipschitz continuous with constants $L_i$, $i=1,...,n$, i.e.
\begin{equation}
|f'_i(x+he_i) - f'_i(x)| \leq L_i |h|, \quad h \in \R, \quad i=1,...,n, \quad x \in Q.
\label{eq:fCoordLipSmGF}
\end{equation}
We set $\beta_i=L_i$, $i=1,...,n$. Then, by definitions in subsection~\ref{S:Not}, we have $\|x\|_E^2 = \sum_{i=1}^n L_i (x^{(i)})^2$, $x \in E$, $d(x) = \frac{1}{2}\sum_{i=1}^n L_i(x^{(i)})^2$, $x \in Q$, $V[z](x)=\frac{1}{2}\sum_{i=1}^n L_i (x^{(i)}-z^{(i)})^2$, $x, z \in Q$. Also, we have $\|g\|_{E,*}^2 = \sum_{i=1}^n L_i^{-1} (g^{(i)})^2$, $g \in E^*$.

We assume that, at any point $x$ in a small vicinity $\bar{Q}$ of the set $Q$, one can calculate an inexact value $\tf(x)$ of the function $f$, s.t. $|\tf(x)-f(x)|\leq \Delta$, $x \in \bar{Q}$. To approximate the gradient of $f$, we use 
$$
\tnf(x) = n\frac{\tf(x+\tau e_i)-\tf(x)}{\tau}e_i, 
$$ 
where $\tau > 0$ is small parameter, which will be chosen later, and the coordinate $i$ is chosen from $i=1,...,n$ randomly with uniform probability $\frac{1}{n}$. Since, we are in the Euclidean setting, we consider $e_i$ also as an element of $E^*$. 

Let us check the assumptions stated in subsection~\ref{S:PrSt&Asmpt}.

\textbf{Randomized Inexact Oracle.}
First, let us show that the finite-difference approximation for the gradient of $f$ can be expressed in the form of \eqref{eq:IODef}. We have
\begin{align*}
\tnf(x) &= n\frac{\tf(x+\tau e_i)-\tf(x)}{\tau}e_i\\
 & = n \left(\langle \nabla f(x), e_i\rangle  + \frac{1}{\tau}( \tf(x+\tau e_i)-\tf(x) - \tau \langle \nabla f(x), e_i\rangle)  \right)e_i.  \notag  
\end{align*}
Taking $\rho = n$, $H = \R$, $\Rf^T: E^* \to \R$ is given by $\Rf^T g = \la g, e_i \ra$, $g \in E^*$, $\Rb: \R \to E^*$ is given by $\Rb t = te_i$, $t \in \R$, we obtain
$$
\tnf(x) = n (\la \nabla f(x), e_i \ra + \xi(x))e_i,
$$
where $\xi(x) = \frac{1}{\tau}( \tf(x+\tau e_i)-\tf(x) - \tau \langle \nabla f(x), e_i\rangle)$.
One can prove that 
\begin{equation*}
	\mE_i n \la \nabla f(x), e_i \ra e_i =  n\E_ie_ie_i^T\nabla f(x) = \nabla f(x), \quad x \in Q,
\end{equation*}
 and, thus, \eqref{eq:IOProp1} holds. It remains to prove \eqref{eq:IOProp2}, i.e., find $\delta$ s.t. for all $x \in Q$, we have $\|\Rb \xi(x)\|_{E,*} \leq \delta$.
\begin{align}
\|\Rb \xi(x)\|_{E,*}  & =  \frac{1}{\sqrt{L_i}}|\xi(x)|  = \frac{1}{\sqrt{L_i}}\left|\frac{1}{\tau}( \tf(x+\tau e_i)-\tf(x) - \tau \langle \nabla f(x), e_i\rangle)  \right|\notag  \\
& = \frac{1}{\tau\sqrt{L_i}} \left| \tf(x+\tau e_i) - f(x+\tau e_i) - (\tf(x) - f(x))\right| \notag \\
& \hspace{1em} + \frac{1}{\tau\sqrt{L_i}} \left|f(x+\tau e_i) - f(x) - \tau \langle \nabla f(x), e_i\rangle \right| \notag \\
& \leq \frac{2\Delta}{\sqrt{L_i}\tau}+ \frac{\tau\sqrt{L_i}}{2}. \notag
\end{align}
Here we used that $|\tf(x)-f(x)|\leq \Delta$, $x \in \bar{Q}$ and \eqref{eq:fCoordLipSm2}, which follows from \eqref{eq:fCoordLipSmGF}. So, we obtain that \eqref{eq:IOProp2} holds with $\delta = \frac{2\Delta}{\sqrt{L_i}\tau}+ \frac{\sqrt{L_i}\tau}{2}$. To balance both terms, we choose $\tau = 2 \sqrt{\frac{\Delta}{L_i}} \leq 2 \sqrt{\frac{\Delta}{L_0}}$, where $L_0 = \min_{i=1,...,n}L_i$.  This leads to equality $\delta = 2 \sqrt{\Delta}$.

\textbf{Regularity of Prox-Mapping.} This assumption can be checked in the same way as in subsection~\ref{S:RCD}.

\textbf{Smoothness.} This assumption can be checked in the same way as in subsection~\ref{S:RCD}.

We have checked that all the assumptions listed in subsection~\ref{S:PrSt&Asmpt} hold. Thus, we can obtain the following convergence rate result for random derivative-free coordinate descent as a corollary of \ag{Theorem}~\ref{Th:1} and \ag{Lemma}~\ref{Lm:AkGrowth}. 
\begin{corollary}
	\label{Cor:UARMRDFCD}
	Let Algorithm~\ref{Alg:UARM} with $\tnf(x) = n\frac{\tf(x+\tau e_i)-\tf(x)}{\tau}e_i$, where $i$ is random and uniformly distributed in $1,...,n$, be applied to Problem \eqref{eq:PrSt} in the setting of this subsection. Let $f_*$ be the optimal objective value and $x_*$ be an optimal point in Problem \eqref{eq:PrSt}. Assume that function value error $\tf(x)-f(x)$ satisfies $|\tf(x)-f(x)| \leq \Delta$, $x \in \bar{Q}$. Denote 
	$$
	P_0^2 = \left(1-\frac{1}{n}\right)(f(x_0)-f_*) + \sum_{i=1}^n\frac{L_i}{2}(u_0^{(i)}-x_*^{(i)})^2.
	$$
	\begin{enumerate}
		\item If the error in the value of the objective $f$ can be controlled and, on each iteration, the error level $\Delta$ satisfies  
		$$
		\Delta \leq \frac{P_0^2}{64n^2A_k^2}, 
		$$
		and $\tau = 2 \sqrt{\frac{\Delta}{L_0}}$ then, for all $k \geq 1$,
		$$
		\mE f(x_{k}) - f_*  \leq \frac{6n^2P_0^2}{(k-1+2n)^2},
		$$
		where $\mE$ denotes the expectation with respect to all the randomness up to step $k$.
		\item If the error in the value of the objective $f$ can not be controlled and $\tau = 2 \sqrt{\frac{\Delta}{L_0}}$, then, for all $k \geq 1$,
		$$
		\mE f(x_{k}) - f_*  \leq \frac{8n^2P_0^2}{(k-1+2n)^2} + 16 (k-1+2n)^2 \Delta. 
		$$
	\end{enumerate}
\end{corollary}

\begin{remark}
	According to \ag{Remark}~\ref{Rm:Rate} and due to the relation $\delta = 2 \sqrt{\Delta}$, we obtain that the error level in the function value should satisfy
	$$
	\Delta \leq \frac{\e^2 }{144 n^2P_0^2}. 
	$$
	The parameter $\tau$ should satisfy
	$$
	\tau \leq \frac{\e }{6 n P_0 \sqrt{L_0}}.
	$$
	At the same time, to obtain an $\e$-solution for Problem \eqref{eq:PrSt}, it is enough to choose
	$$
	k = \max\left\{\left\lceil n\sqrt{\frac{6P_0^2}{\e}} + 1 -2 n \right\rceil,0\right\}.
	$$
\end{remark}

\subsection{Accelerated Random Derivative-Free Block-Coordinate Descent}
\label{S:RDFBCD}
In this subsection, we consider the same setting as in subsection~\ref{S:RBCD}, except for Randomized Inexact Oracle. Instead of block derivative, we use here its finite-difference approximation. 
As in subsection~\ref{S:RBCD}, we consider Euclidean setup and entropy setup. 

\textit{Euclidean setup.} We assume that, for all $i=1,...,n$, $E_i=\R^{p_i}$; $Q_i$ is a simple closed convex set; $\|x^{(i)}\|_i^2=\la B_i x^{(i)}, x^{(i)} \ra$, $x^{(i)} \in E_i$, where $B_i$ is symmetric positive semidefinite matrix; $d_i(x^{(i)})=\frac12\|x^{(i)}\|_i^2$, $x^{(i)} \in Q_i$, and, hence, $V_i[z^{(i)}](x^{(i)}) = \frac12\|x^{(i)}-z^{(i)}\|_i^2$, $x^{(i)}, z^{(i)} \in Q_i$.

\textit{Entropy setup.} We assume that, for all $i=1,...,n$, $E_i=\R^{p_i}$; $Q_i$ is standard simplex in $\R^{p_i}$, i.e., $Q_i = \{x^{(i)} \in \R^{p_i}_+: \sum_{j=1}^{p_i}[x^{(i)}]_j = 1\}$; $\|x^{(i)}\|_i=\|x^{(i)}\|_1=\sum_{j=1}^{p_i}|[x^{(i)}]_j|$, $x^{(i)} \in E_i$; $d_i(x^{(i)})=\sum_{j=1}^{p_i}[x^{(i)}]_j \ln [x^{(i)}]_j$, $x^{(i)} \in Q_i$, and, hence, $V_i[z^{(i)}](x^{(i)}) = \sum_{j=1}^{p_i}[x^{(i)}]_j \ln \frac{[x^{(i)}]_j}{[z^{(i)}]_j}$, $x^{(i)}, z^{(i)} \in Q_i$. 

Note that, in each block, one also can choose other proximal setups from \cite{ben-tal2015lectures}. Combination of different setups in different blocks is also possible, i.e., in one block it is possible to choose the Euclidean setup and in an another block one can choose the entropy setup.

Using operators $U_i$, $i=1,...,n$ defined in \eqref{eq:UDef}, for each $i=1,...,n$, the $i$-th block derivative of $f$ can be written as $f_i'(x) = U_i^T \nabla f(x)$.
We assume that the gradient of $f$ in \eqref{eq:PrSt} is block-wise Lipschitz continuous with constants $L_i$, $i=1,...,n$ with respect to chosen norms $\|\cdot \|_i$, i.e.,
\begin{equation}
\|f'_i(x+U_ih^{(i)}) - f'_i(x)\|_{i,*} \leq L_i \|h^{(i)}\|_i, \quad h^{(i)} \in E_i, \quad i=1,...,n \quad x \in Q.
\label{eq:fBlockLipSmGF}
\end{equation}
We set $\beta_i=L_i$, $i=1,...,n$. Then, by definitions in subsection~\ref{S:Not}, we have $\|x\|_E^2 = \sum_{i=1}^n L_i \|x^{(i)}\|_2^2$, $x \in E$, $d(x) = \frac12 \sum_{i=1}^n L_i\|x^{(i)}\|_2^2$, $x \in Q$, $V[z](x)=\frac12\sum_{i=1}^n L_i \|x^{(i)}-z^{(i)} \|_2^2$, $x, z \in Q$. Also, we have $\|g\|_{E,*}^2 = \sum_{i=1}^n L_i^{-1} \|g^{(i)}\|_2^2$, $g \in E^*$.

We assume that, at any point $x$ in a small vicinity $\bar{Q}$ of the set $Q$, one can calculate an inexact value $\tf(x)$ of the function $f$, s.t. $|\tf(x)-f(x)|\leq \Delta$, $x \in \bar{Q}$. To approximate the gradient of $f$, we use 
\begin{align}
\tnf(x) = n\widetilde{U}_i \left( \frac{\tf(x+\tau U_ie_1)-\tf(x)}{\tau},...,\frac{\tf(x+\tau U_ie_{p_i})-\tf(x)}{\tau}\right)^T, 
\label{eq:RDFBCDtnfDef}  
\end{align}
where $\tau > 0$ is small parameter, which will be chosen later, a block number $i$ is chosen from $i=1,...,n$ randomly with uniform probability $\frac{1}{n}$, $e_1,...,e_{p_i}$ are coordinate vectors in $E_i$, $U_i$ is defined in \eqref{eq:UDef}, $\widetilde{U}_i$ is defined in \eqref{eq:tUDef}.

Let us check the assumptions stated in subsection~\ref{S:PrSt&Asmpt}.

\textbf{Randomized Inexact Oracle.}
First, let us show that the random derivative-free block-coordinate approximation for the gradient of $f$ can be expressed in the form of \eqref{eq:IODef}. Denote $\tilde{g}_i = \frac{1}{\tau} \left( \tf(x+\tau U_i e_1)-\tf(x),...,\tf(x+\tau U_i e_{p_i})-\tf(x)\right)^T \in E_i$, $i=1,...,n$. We have
\begin{align}
\tnf(x) & = n\widetilde{U}_i\tilde{g}_i  = n\widetilde{U}_i \left( U_i^T\nabla f(x) + (\tilde{g}_i - U_i^T\nabla f(x))  \right).  
\notag
\end{align}
Taking $\rho = n$, $H = E_i$, $\Rf^T: E^* \to E_i^*$ be given by $\Rf^T g = U_i^T g$, $g \in E^*$ and $\Rb: E_i^* \to E^*$ be given by $\Rb g^{(i)} = \widetilde{U}_i g^{(i)}$, $g^{(i)} \in E_i^*$, we obtain
$$
\tnf(x) = n \widetilde{U}_i(U_i^T\nabla f(x) + \xi(x)),
$$
where $\xi(x) = \tilde{g}_i - U_i^T\nabla f(x)$. Since $i \in R[1,n]$, one can prove that $\mE_i n \widetilde{U}_iU_i^T \nabla f(x) =  \nabla f(x)$, $x \in Q$, and, thus, \eqref{eq:IOProp1} holds. 
It remains to prove \eqref{eq:IOProp2}, i.e., find $\delta$ s.t. for all $x \in Q$, we have $\|\Rb \xi(x)\|_{E,*} \leq \delta$. 
Let us fix any $i$ from $1,...,n$ and $j$ from $1,...,p_i$. Then, for any $x \in \bar{Q}$, the $j$-th coordinate of $\xi(x) = \tilde{g}_i - U_i^T\nabla f(x)$ can be estimated as follows
\begin{align}
|[\xi(x)]_j| &= \left|\frac{1}{\tau}( \tf(x+\tau U_i e_j)-\tf(x) - \tau \la U_i^T \nabla f(x), e_j\ra \right| \notag \\
& = \frac{1}{\tau }\left| \tf(x+\tau U_i e_j) - f(x+\tau U_i e_j) - (\tf(x) - f(x))\right| \notag \\
& \hspace{1em} + \frac{1}{\tau }\left| (f(x+\tau U_i e_j) - f(x) - \tau \la U_i^T \nabla f(x), e_j\ra ) \right| \notag \\
& \leq \frac{2\Delta }{\tau }+ \frac{\tau L_i}{2}.
\label{eq:xijEst}
\end{align}
Here we used that $|\tf(x)-f(x)|\leq \Delta$, $x \in \bar{Q}$, \eqref{eq:fBlockLipSm2}, which follows from \eqref{eq:fBlockLipSmGF}. 
In our setting, for any $i = 1,...,n$, $\|\cdot\|_{i,*}$ is either $\max$-norm (for the entropy case) or Euclidean norm (for the Euclidean case). Thus, in the worst case of Euclidean norm
\begin{align*}
 \|\Rb \xi(x)\|_{E,*}   &=  \|\widetilde{U}_i \xi(x)\|_{E,*}\\
&  = \frac{1}{\sqrt{L_i}} \left\|\xi(x) \right\|_{i,*}\\ &\stackrel{\eqref{eq:xijEst}}{\leq} \frac{\sqrt{p_i}}{\sqrt{L_i}} \left(\frac{2\Delta }{\tau }+ \frac{\tau L_i}{2} \right)\\
& \leq  \sqrt{p_{\max}} \left(\frac{2\Delta }{\tau \sqrt{L_i}}+ \frac{\tau \sqrt{L_i}}{2} \right), \notag 
\end{align*}
where $p_{\max} = \max_{i=1,...,n} p_i$.
So, we obtain that \eqref{eq:IOProp2} holds with 
\begin{align*}
	\delta = \sqrt{p_{\max}} \left(\frac{2\Delta }{\tau \sqrt{L_i}}+ \frac{\tau \sqrt{L_i}}{2} \right).
\end{align*}
 To balance both terms we choose $\tau =  2 \sqrt{\frac{\Delta}{L_i}} \leq 2 \sqrt{\frac{\Delta}{L_0}}$, where $L_0 = \min_{i=1,...,n}L_i$. This leads to equality $\delta = 2  \sqrt{p_{\max} \Delta}$.

\textbf{Regularity of Prox-Mapping.} This assumption can be checked in the same way as in subsection~\ref{S:RBCD}.

\textbf{Smoothness.} This assumption can be checked in the same way as in subsection~\ref{S:RBCD}.

We have checked that all the assumptions listed in subsection~\ref{S:PrSt&Asmpt} hold. Thus, we can obtain the following convergence rate result for random derivative-free block-coordinate descent as a corollary of \ag{Theorem}~\ref{Th:1} and \ag{Lemma}~\ref{Lm:AkGrowth}. 
\begin{corollary}
	\label{Cor:UARMBlockCoordGF}
	Let Algorithm~\ref{Alg:UARM} with $\tnf(x)$ defined in \eqref{eq:RDFBCDtnfDef}, be applied to Problem ~\eqref{eq:PrSt} in the setting of this subsection. Let $f_*$ be the optimal objective value and $x_*$ be an optimal point in Problem ~\eqref{eq:PrSt}. Assume that function value error $\tf(x)-f(x)$ satisfies $|\tf(x)-f(x)| \leq \Delta$, $x \in \bar{Q}$. Denote 
	$$
	P_0^2 = \left(1-\frac{1}{n}\right)(f(x_0)-f_*) + V[u_0](x_*).
	$$
	\begin{enumerate}
		\item If the error in the value of the objective $f$ can be controlled and, on each iteration, the error level $\Delta$ satisfies 
		$$
		\Delta \leq \frac{P_0^2}{64n^2p_{\max}A_k^2} , 
		$$
		and $\tau = 2 \sqrt{\frac{\Delta}{L_0}}$ then, for all $k \geq 1$,
		$$
		\mE f(x_{k}) - f_*  \leq \frac{6n^2P_0^2}{(k-1+2n)^2} ,
		$$
		where $\mE$ denotes the expectation with respect to all the randomness up to step $k$.
		\item If the error in the value of the objective $f$ can not be controlled and $\tau = 2 \sqrt{\frac{\Delta}{L_0}}$, then, for all $k \geq 1$,
		$$
		\mE f(x_{k}) - f_*  \leq \frac{8n^2P_0^2}{(k-1+2n)^2} + 16 (k-1+2n)^2 p_{\max} \Delta. 
		$$
	\end{enumerate}
\end{corollary}

\begin{remark}
	According to \ag{Remark}~\ref{Rm:Rate} and due to the relation $\delta = 2 \sqrt{p_{\max}\Delta}$, we obtain that the error level in the function value should satisfy
	$$
	\Delta \leq \frac{\e^2 }{144 n^2p_{\max}P_0^2}. 
	$$
	The parameter $\tau$ should satisfy
	$$
	\tau \leq \frac{\e }{6 n P_0 \sqrt{L_0}}.
	$$
	At the same time, to obtain an $\e$-solution for Problem ~\eqref{eq:PrSt}, it is enough to choose
	$$
	k = \max\left\{\left\lceil n\sqrt{\frac{6P_0^2}{\e}} + 1 -2 n \right\rceil,0\right\}.
	$$
\end{remark}

\subsection{Accelerated Random Derivative-Free Block-Coordinate Descent with Random Approximations for Block Derivatives}
\label{S:RDFBCDRA}
In this subsection, we combine random block-coordinate descent of subsection~\ref{S:RBCD} with random derivative-free directional search described in subsection~\ref{S:RDFDS} and random derivative-free coordinate descent described in \ref{S:RDFCD}. 
We construct randomized approximations for block derivatives based on finite-difference approximation of directional derivatives.
Unlike subsection~\ref{S:RBCD}, we consider only Euclidean setup. We assume that, for all $i=1,...,n$, $E_i=\R^{p_i}$; $\|x^{(i)}\|_i^2=\|x^{(i)}\|_2^2$, $x^{(i)} \in E_i$;  $Q_i$ is either $E_i$, or $\otimes_{j=1}^{p_i}Q_{ij}$, where $Q_{ij} \subseteq \R$ are closed convex sets; $d_i(x^{(i)})=\frac12\|x^{(i)}\|_i^2$, $x^{(i)} \in Q_i$ and, hence, $V_i[z^{(i)}](x^{(i)}) = \frac12\|x^{(i)}-z^{(i)}\|_i^2$, , $x^{(i)}, z^{(i)} \in Q_i$. For the case, $Q_i = E_i$, we consider randomization on the Euclidean sphere of radius 1, as in subsection~\ref{S:RDFDS}. For the case, $Q_i = \otimes_{j=1}^{p_i}Q_{ij}$, we consider coordinate-wise randomization, as in subsection~\ref{S:RDFCD}. 

Using operators $U_i$, $i=1,...,n$ defined in ~\eqref{eq:UDef}, for each $i=1,...,n$, the $i$-th block derivative of $f$ can be written as $f_i'(x) = U_i^T \nabla f(x)$.
We assume that the gradient of $f$ in ~\eqref{eq:PrSt} is block-wise Lipschitz continuous with constants $L_i$, $i=1,...,n$ with respect to chosen norms $\|\cdot \|_i$, i.e.,
\begin{equation}
\|f'_i(x+U_ih^{(i)}) - f'_i(x)\|_{i,*} \leq L_i \|h^{(i)}\|_i, \quad h^{(i)} \in E_i, \quad i=1,...,n \quad x \in Q.
\label{eq:fBlockLipSmGFRA}
\end{equation}
We set $\beta_i=L_i$, $i=1,...,n$. Then, by definitions in subsection~\ref{S:Not}, we have $\|x\|_E^2 = \sum_{i=1}^n L_i \|x^{(i)}\|_2^2$, $x \in E$, $d(x) = \frac12 \sum_{i=1}^n L_i\|x^{(i)}\|_2^2$, $x \in Q$, $V[z](x)=\frac{1}{2}\sum_{i=1}^n L_i \|x^{(i)}-z^{(i)} \|_2^2$, $x, z \in Q$. Also, we have $\|g\|_{E,*}^2 = \sum_{i=1}^n L_i^{-1} \|g^{(i)}\|_2^2$, $g \in E^*$.

We assume that, at any point $x$ in a small vicinity $\bar{Q}$ of the set $Q$, one can calculate an inexact value $\tf(x)$ of the function $f$, s.t. $|\tf(x)-f(x)|\leq \Delta$, $x \in \bar{Q}$. To approximate the gradient of $f$, we first randomly choose a block $i \in 1,...,n$ with probability $p_i/p$, where $p=\sum_{i=1}^np_i$. Then we use one the following types of random directions $e \in E_i$ to approximate the block derivative $f'_i(x)$ by a finite-difference. 
\begin{enumerate}
	\item If $Q_i = E_i$, we take $e \in E_i$ to be random vector uniformly distributed on the Euclidean sphere of radius 1, i.e. $\Sp_2(1):=\{s \in \R^{p_i}: \|s\|_2=1\}$. We call this \textit{unconstrained case}.
	\item If $Q_i = \otimes_{j=1}^{p_i}Q_{ij}$, we take $e$ to be random uniformly chosen from $1,...,p_i$ coordinate vector, i.e. $e = e_j \in E_i$ with probability $\frac{1}{p_i}$. We call this \textit{separable case}.
\end{enumerate}
Based on these randomizations and inexact function values, our randomized approximation for the gradient of $f$ is 
$$
\tnf(x) = p \widetilde{U}_i \frac{\tf(x+\tau U_i e)-\tf(x)}{\tau}e, 
$$ 
where $\tau > 0$ is small parameter, which will be chosen later, $U_i$ is defined in ~\eqref{eq:UDef} and $\widetilde{U}_i$ is defined in ~\eqref{eq:tUDef}.

Let us check the assumptions stated in subsection~\ref{S:PrSt&Asmpt}.

\textbf{Randomized Inexact Oracle.}
First, let us show that the random derivative-free block-coordinate approximation for the gradient of $f$ can be expressed in the form of ~\eqref{eq:IODef}. We have
\begin{align}
\tnf(x) & = p\widetilde{U}_i\frac{\tf(x+\tau U_i e)-\tf(x)}{\tau}e \notag \\
& = p\widetilde{U}_i \left( \langle U_i^T\nabla f(x), e\rangle e + \frac{1}{\tau}( \tf(x+\tau U_i e)-\tf(x) - \tau \langle U_i^T \nabla f(x), e\rangle) e \right).  \label{eq:RDFBCDRAtnfDef}  
\end{align}
Taking $\rho = p$, $H = E_i$, $\Rf^T: E^* \to E_i^*$ be given by $\Rf^T g = \la U_i^T g,e \ra e$, $g \in E^*$ and $\Rb: E_i^* \to E^*$ be given by $\Rb g^{(i)} = \widetilde{U}_i g^{(i)}$, $g^{(i)} \in E_i^*$, we obtain
$$
\tnf(x) = p \widetilde{U}_i(\langle U_i^T\nabla f(x), e\rangle e + \xi(x)),
$$
where $\xi(x) = \frac{1}{\tau}( \tf(x+\tau U_i e)-\tf(x) - \tau \langle U_i^T \nabla f(x), e\rangle) e$. By the choice of probability distributions for $i$ and $e$ and their independence, we have, for all $x \in Q$, 
\begin{align}
\E_{i,e} p \widetilde{U}_i\langle U_i^T\nabla f(x), e\rangle e & =  p\E_{i,e}  \widetilde{U}_ie e^T U_i^T\nabla f(x) = p\E_{i}  \widetilde{U}_i \left(E_{e}e e^T\right) U_i^T\nabla f(x) \notag \\
&= p\E_{i} \frac{1}{p_i} \widetilde{U}_i U_i^T \nabla f(x) = \nabla f(x) \notag 
\end{align}
and, thus, ~\eqref{eq:IOProp1} holds. 
It remains to prove ~\eqref{eq:IOProp2}, i.e., find $\delta$ s.t. for all $x \in Q$, we have $\|\Rb \xi(x)\|_{E,*} \leq \delta$. We have
\begin{align*}
& \|\Rb \xi(x)\|_{E,*}   =  \|\widetilde{U}_i \xi(x)\|_{E,*}  = \frac{1}{\sqrt{L_i}} \left\|\frac{1}{\tau}( \tf(x+\tau U_i e)-\tf(x) - \tau \langle U_i^T \nabla f(x), e\rangle) e \right\|_{i,*}\notag  \\
& =\frac{\left\|e \right\|_{i,*} }{\tau \sqrt{L_i}}\Big| \tf(x+\tau U_i e) - f(x+\tau U_i e) - (\tf(x) - f(x))\\
&\hspace{1em} + (f(x+\tau U_i e) - f(x) - \tau \langle U_i^T \nabla f(x), e\rangle ) \Big| \notag \\
& \leq \frac{2\Delta \left\|e \right\|_{i,*}}{\tau \sqrt{L_i}}+ \frac{\tau  \left\|e \right\|_{i,*} \left\|e \right\|_{i}^2 \sqrt{L_i}}{2} 
\notag \\
& = \frac{2\Delta }{\tau \sqrt{L_i}}+ \frac{\tau \sqrt{L_i}}{2}. \notag
\end{align*}
Here we used that $|\tf(x)-f(x)|\leq \Delta$, $x \in \bar{Q}$, ~\eqref{eq:fBlockLipSm2}, which follows from ~\eqref{eq:fBlockLipSmGFRA}, and that the norms $\|\cdot\|_i$, $\|\cdot\|_{i,*}$ are standard Euclidean. So, we obtain that ~\eqref{eq:IOProp2} holds with $\delta = \frac{2\Delta }{\tau \sqrt{L_i}}+ \frac{\tau \sqrt{L_i}}{2}$. To balance both terms we choose $\tau =  2 \sqrt{\frac{\Delta}{L_i}} \leq 2 \sqrt{\frac{\Delta}{L_0}}$, where $L_0 = \min_{i=1,...,n}L_i$. This leads to equality $\delta = 2  \sqrt{\Delta}$.

\textbf{Regularity of Prox-Mapping.} Separable structure of $Q$ and $V[u](x)$ means that the problem ~\eqref{eq:uPlusDef} boils down to $n$ independent problems of the form 
$$
u_+^{(l)} = \arg \min_{x^{(l)} \in Q_l} \left\{\frac{L_l}{2}\|u^{(l)}-x^{(l)}\|_2^2+ \alpha \la U_l^T\tnf(y),x^{(l)} \ra\right\}, \quad l = 1,...,n.
$$
Since $\tnf(y)$ has non-zero components only in the block $i$, $U_l^T\tnf(y)$ is zero for all $l \ne i$. Thus, $u-u_+$ has non-zero components only in the block $i$ and $U_i(u^{(i)}-u_+^{(i)}) = u - u_+$. 

In the unconstrained case, similarly to subsection~\ref{S:RDS}, we obtain that $u^{(i)}-u_+^{(i)}=\gamma e$, where $\gamma$ is some constant. Using these two facts, we obtain
\begin{align}
\la \Rb\Rf^T\nabla f(y), u - u_+ \ra & = \la \widetilde{U}_i\la U_i^T\nabla f(y), e \ra e , u - u_+ \ra \notag \\
& = \la\la U_i^T\nabla f(y), e \ra e ,  \widetilde{U}_i^T(u - u_+) \ra \notag \\
& = \la\la U_i^T\nabla f(y), e \ra e ,  u^{(i)}-u_+^{(i)} \ra \notag \\
& = \la\la U_i^T\nabla f(y), e \ra e ,  \gamma e \ra \notag \\
& = \la U_i^T\nabla f(y), \gamma e \ra \la e ,   e \ra \notag \\
& = \la U_i^T\nabla f(y), u^{(i)}-u_+^{(i)}  \ra \notag \\
& = \la \nabla f(y), U_i(u^{(i)}-u_+^{(i)}) \ra \notag \\
& = \la \nabla f(y), u - u_+ \ra,  \notag
\end{align}
which proves ~\eqref{eq:RPMA} for the unconstrained case.

In the separable case, similarly to subsection~\ref{S:RCD}, we obtain that $u^{(i)}-u_+^{(i)}$ has only one $j$-th non-zero coordinate, where $j \in 1,...,p_i$. Hence, $\la e_j, u^{(i)}-u_+^{(i)} \ra e_j = u^{(i)}-u_+^{(i)}$. So, we get, 
\begin{align}
\la \Rb\Rf^T\nabla f(y), u - u_+ \ra & = \la \widetilde{U}_i\la U_i^T\nabla f(y), e_j \ra e_j , u - u_+ \ra \notag \\
& = \la\la U_i^T\nabla f(y), e_j \ra e_j ,  \widetilde{U}_i^T(u - u_+) \ra \notag \\
& = \la\la U_i^T\nabla f(y), e_j \ra e_j ,  u^{(i)}-u_+^{(i)} \ra \notag \\
& = \la U_i^T\nabla f(y), e_j \ra \la e_j ,  u^{(i)}-u_+^{(i)} \ra \notag \\
& = \la U_i^T\nabla f(y), \la e_j ,  u^{(i)}-u_+^{(i)} \ra e_j \ra  \notag \\
& = \la U_i^T\nabla f(y), u^{(i)}-u_+^{(i)}  \ra \notag \\
& = \la \nabla f(y), U_i(u^{(i)}-u_+^{(i)}) \ra \notag \\
& = \la \nabla f(y), u - u_+ \ra,  \notag
\end{align}
which proves ~\eqref{eq:RPMA} for the separable case.

\textbf{Smoothness.} This assumption can be checked in the same way as in subsection~\ref{S:RBCD}.

We have checked that all the assumptions listed in subsection~\ref{S:PrSt&Asmpt} hold. Thus, we can obtain the following convergence rate result for random derivative-free block-coordinate descent with random approximations for block derivatives as a corollary of \ag{Theorem}~\ref{Th:1} and \ag{Lemma}~\ref{Lm:AkGrowth}. 
\begin{corollary}
	\label{Cor:UARMBlockCoordGFRA}
	Let \ref{Alg:UARM} with $\tnf(x)$ defined in ~\eqref{eq:RDFBCDRAtnfDef}, be applied to Problem ~\eqref{eq:PrSt} in the setting of this subsection. Let $f_*$ be the optimal objective value and $x_*$ be an optimal point in Problem ~\eqref{eq:PrSt}. Assume that function value error $\tf(x)-f(x)$ satisfies $|\tf(x)-f(x)| \leq \Delta$, $x \in \bar{Q}$. Denote 
	$$
	P_0^2 = \left(1-\frac{1}{p}\right)(f(x_0)-f_*) + \sum_{i=1}^n\frac{L_i}{2}\|u_0^{(i)}-x_*^{(i)}\|_2^2.
	$$
	\begin{enumerate}
		\item If the error in the value of the objective $f$ can be controlled and, on each iteration, the error level $\Delta$ satisfies 
		$$
		\Delta \leq \frac{P_0^2}{64p^2A_k^2} , 
		$$
		and $\tau = 2 \sqrt{\frac{\Delta}{L_0}}$ then, for all $k \geq 1$,
		$$
		\mE f(x_{k}) - f_*  \leq \frac{6p^2P_0^2}{(k-1+2p)^2} ,
		$$
		where $\mE$ denotes the expectation with respect to all the randomness up to step $k$.
		\item If the error in the value of the objective $f$ can not be controlled and $\tau = 2 \sqrt{\frac{\Delta}{L_0}}$, then, for all $k \geq 1$,
		$$
		\mE f(x_{k}) - f_*  \leq \frac{8p^2P_0^2}{(k-1+2p)^2} + 16 (k-1+2p)^2 \Delta. 
		$$
	\end{enumerate}
\end{corollary}

\begin{remark}
	According to \ag{Remark}~\ref{Rm:Rate} and due to the relation $\delta = 2 \sqrt{\Delta}$, we obtain that the error level in the function value should satisfy
	$$
	\Delta \leq \frac{\e^2 }{144 p^2P_0^2}. 
	$$
	The parameter $\tau$ should satisfy
	$$
	\tau \leq \frac{\e }{6 p P_0 \sqrt{L_0}}.
	$$
	At the same time, to obtain an $\e$-solution for Problem ~\eqref{eq:PrSt}, it is enough to choose
	$$
	k = \max\left\{\left\lceil p\sqrt{\frac{6P_0^2}{\e}} + 1 -2 p \right\rceil,0\right\}.
	$$
\end{remark}

\section{Model generality in a Non-accelerated Random Block-Coordinate Descent}
\ag{Above we propose a unified view on different randomized gradient-free and (block) coordinate-wise schemes. In the last decade an interest to structural optimization has grown \cite{nesterov2018lectures}. The most general results (see  \cite{dvinskikh2020accelerated,gasnikov2017universal,gasnikov2019fast,stonyakin2019inexact,stonyakin2020inexact,stonyakin2019gradient}) allow to unify different structures of the problem (composite problems, $\max$-type problems e.t.c.) in one envelop (model generality \cite{gasnikov2017universal}). But this envelop deals with the first-order methods. Motivated by Random Block Coordinate Descent (RBCD)  \cite{nesterov2017random} we try to develop RBCD in model generality. We expect that further unification of model generality conception \cite{gasnikov2017universal} with the results of this paper is also possible.}

\ag{Typical randomized algorithm has such kind of auxiliary problem
$$x_{k+1}=\arg \min_{x \in Q} \{V[x_{k}](x) + \alpha \la \tnf(x_{k}), x - x_k \ra \},$$
see for example \eqref{eq:uPlusDef} and \eqref{eq:ukp1Def} in Algorithm~\ref{Alg:UARM}. Here the term $\la \tnf(x_{k}), x - x_k \ra$ is a (random) linear <<model>> of target function in considered point $x_k$.
The idea is to replace this model $\la \tnf(x_{k}), x - x_k \ra$ on something more general. 
For example, if $f(x):= f(x) + h(x)$ we may expect, that composite model $\la \tnf(x_{k}), x - x_k \ra + h(x) - h(x_k)$ is
also well suited and may lead to a better results than standard model $\la \tnfg(x_{k}), x - x_k \ra$.}


\ag{Assume that}


\begin{equation*}\label{eq:model_generality}
    \ag{V[z](x)=\frac{1}{2}\sum_{i=1}^n L_i \|x^{(i)}-z^{(i)} \|_2^2, \quad Q = \otimes_{i=1}^{n}Q_{i}}.
\end{equation*}

\ag{Function $\psi_i(y,x)$ determines $i$-th part of the model of target function, that corresponds to block $Q^{(i)}$. For example, in composite case if $h(x) = \sum_{i=1}^n h_i(x^{(i)})$ we can consider $$\psi_i(y,x) = \la \widetilde{U_i} U_i^T \nabla f(x), y - x \ra + h_i(y^{(i)}) - h_i(x^{(i)}).$$}

\begin{assumption}[convexity]\label{A:MGEM}
   Assume \ag{that for all $x \in Q$ model $\psi_i(y, x)$ is convex function of $y\in Q$ and there exists such $\gamma > 0$} that for any $x, y \in Q$
\begin{equation}\label{eq:mgfp}
    \mE_i \psi_i(y, x) \leq \frac{1}{\gamma}(f(y) - f(x)).
\end{equation} 
\end{assumption}

Parameter $\gamma = n$ for composite case.

\begin{assumption}[smoothness]\label{A:MGLI}
   \ag{F}or any 
   $\ag{x_k,x_{k+1}}$ \ag{generated by Algorithm~\ref{Alg:NRBCD}} 
\begin{equation}\label{eq:ng_pi}
     f(x_{k+1}) \leq f(x_k) + \psi_i(x_{k+1}, x_k) + V[x_k](x_{k+1}).
\end{equation}
\end{assumption}
\ag{Assumption~\ref{A:MGLI} typically holds true with additional condition that $x_k,x_{k+1}$ generated by Algorithm~\ref{Alg:NRBCD} satisfy   $\ag{x_k^{(j)} = x_{k+1}^{(j)}}$ for $j \neq i$.} 
\begin{algorithm}[h!]
\caption{Non-accelerated Random Block-Coordinate Descent (NRBCD)}
\label{Alg:NRBCD}

 \KwIn{starting point $x_0 \in Q^0 = \otimes_{i=1}^n Q_i^0$, \at{number of iterations $N$,} prox-setup: $d(x)$, $V[u] (x)$, see subsection~\ref{S:Not}.}

 Set $k=0$.
 
 \Repeat{\ag{$k\le N$}}{Choose randomly $i \in \{1, \ldots, n\}$ \ag{($\mathds{P}(i = j) = 1/n$ for all $j=1,...,n$)} \\
 \begin{equation}\label{eq:NRBCD_prox}
				x_{k+1} = \ag{\arg \min_{x \in Q}\{V[x_k](x) + \psi_i (x, x_k)\}}
				\end{equation}
 Set $k=k+1$.}
 
 \KwOut{The point $x_{N}$.}

\end{algorithm}


\begin{theorem}
    Let the assumptions~\ref{A:MGEM} and~\ref{A:MGLI} hold. Let $x_k$ be generated by Algorithm~\ref{Alg:NRBCD}. Let $f_*$ be the optimal objective value and $x^*$ be an optimal point in Problem~\eqref{eq:PrSt}. Then, for all $N \geq 1$,
    \begin{equation*}
        \mE f(\overline{x}_N) - f(x^*) \leq \frac{\ag{\gamma}}{N}(f(x_0) - f(x^*)) + \frac{\ag{\gamma}}{N}V[x_0](x^*),
    \end{equation*}
    where $\overline{x}_N = \frac{1}{N}\sum_{k=0}^{N-1} x_k$.
\end{theorem}
\begin{proof}
    
\ag{From \eqref{eq:ng_pi}, \eqref{eq:NRBCD_prox} follows formula (3.12)  of \cite{gasnikov2017universal} 
\begin{equation*}
-\psi_i(x, x_k) \leq  f(x_k) - f(x_{k+1}) + V[x_k](x) - V[x_{k+1}](x)
\end{equation*}
for all $x\in Q$.  From this formula and \eqref{eq:mgfp}
\begin{align*}
\frac{1}{\gamma}(f(x_k) - f(x)) \leq
-\mE_i\psi(x, x_k)  \leq  f(x_k) -\mE_if(x_{k+1}) + V[x_k](x) - \mE_iV[x_{k+1}](x). 
\end{align*}
}

\ag{Then take the full mathematical expectation from each inequality, where} $k = 0, \ldots, \ag{N-1}$, \ag{put} $x = \ag{x_*}$ and sum all the inequalities for $k = 0, \ldots, \ag{N-1}$, \ag{w}e obtain that 
\begin{equation*}
    \frac{1}{\ag{\gamma}}\sum_{k=0}^{\ag{N-1}}\mE f(x_k) - \frac{N}{\ag{\gamma}}f(\ag{x_*}) \leq \mE(f(x_0) - f(x_{\ag{N}})) + V[x_0](\ag{x_*}) - \mE V[x_{\ag{N}}](\ag{x_*})
\end{equation*}
\ag{Since} $\mE V[x_{\ag{N}}](\ag{x_*}) \geq 0$, then
\begin{equation*}
    \frac{1}{\ag{\gamma}}\sum_{k=0}^{\ag{N-1}}\mE f(x_k) - \frac{N}{\ag{\gamma}}f(\ag{x_*}) \leq \mE(f(x_0) - f(x_{\ag{N}})) + V[x_0](\ag{x_*})
\end{equation*}
Using Yensen inequality and the fact that $f(x_{\ag{N}}) \geq f(\ag{x_*})$ we obtain final result
\begin{equation*}
    \mE f(\overline{x}_N) - f(\ag{x_*}) \leq \frac{\ag{\gamma}}{N}(f(x_0) - f(\ag{x_*})) + \frac{\ag{\gamma}}{N}V[x_0](\ag{x_*}).
\end{equation*}
\end{proof}

\section*{Conclusion}
In this paper, we introduce a unifying framework, which allows to construct different types of accelerated randomized methods for smooth convex optimization problems and to prove convergence rate theorems for these methods. As we show, our framework is rather flexible and allows to reproduce known results as well as obtain new methods with convergence rate analysis. At the moment randomized methods for empirical risk minimization problems are not directly covered by our framework. It seems to be a n interesting direction for further research. Another directions, in which we actually work, include generalization of our framework for strongly convex problems based on well-known restart technique. Another direction of our work is connected to non-uniform probabilities for sampling of coordinate blocks and composite optimization problems.

\textbf{Acknowledgments.} The authors are very grateful to Yu. Nesterov and V. Spokoiny for fruitful discussions. \ag{Our interest to this field was initiated by the paper \cite{nesterov2017random}.}

\bibliographystyle{siamplain}
\bibliography{references}


\end{document}